\documentclass{article}	
\usepackage[utf8]{inputenc}				%(only for the pdftex engine)
\usepackage{lmodern,anysize} 
\usepackage{microtype,amsthm,amsfonts}
\usepackage[numbers,square,sort&compress]{natbib}

\usepackage{subcaption}  % Include this in your preamble

\usepackage{graphicx,hyperref}
\usepackage{ltablex}  % Combines longtable and tabularx

\keepXColumns  % Ensure the 'X' columns maintain their behavior across page breaks

\graphicspath{{../.}{logos}}

\usepackage[affil-it]{authblk} % Package to handle multiple authors and affiliations

\begin{document}

%%%%--------------------------------------------%%%
%	\articletype{Research Article}
%	\received{Month	DD, YYYY}
%	\revised{Month	DD, YYYY}
%  \accepted{Month	DD, YYYY}
%  \journalname{De~Gruyter~Journal}
%  \journalyear{YYYY}
%  \journalvolume{XX}
%  \journalissue{X}
%  \startpage{1}
%  \aop
%  \DOI{10.1515/sample-YYYY-XXXX}
%%%%--------------------------------------------%%%

\title{Golf Strategy Optimization\\ and the ``Drive for show, putt for dough'' adage}
\author[1]{Gautier Stauffer}
\author[2]{Matthieu Guillot}
\affil[1]{Faculty of Business and Economics (HEC Lausanne), Department of Operations, University of Lausanne, Quartier Unil-Chamberonne, 1015 Lausanne, Switzerland, Email: \texttt{gautier.stauffer@unil.ch}}
\affil[2]{Laboratoire DISP, IUT Lumière, Université Lumière Lyon2, France, Email: \texttt{matthieu.guillot@univ-lyon2.fr}}

\maketitle

\abstract{%

This study explores strategic decision-making in professional golf’s Stroke Play format through a computational lens. We develop a Markov Decision Process (MDP) model—specifically, a stochastic shortest path formulation—to optimize a golfer’s strategy on any given course, incorporating both course layout and player skill data. While MDPs have been widely used in sports analytics, applying them to golf presents significant scalability challenges due to the curse of dimensionality. Our primary objective is not to predict player performance with high precision, but rather to demonstrate that an exact, data-driven MDP approach is computationally tractable on full scale, real-world instances. We show that, with careful problem structuring, low-level coding, and efficient memory management, it is possible to solve such large-scale models without resorting to heuristics or Q-learning approximations, as used in existing approaches. To illustrate the model’s potential, we show how one can use PGA Tour data and aerial course imagery to simulate strategic outcomes and analyze how different skill profiles influence performance. In particular, we assess the relative impact of driving and putting, challenging the popular adage ``Drive for show, putt for dough.''  These results support the value of our methodology as a robust proof of concept and a foundation for future enhancements. All code and analyses (in R and C++) are made available as open-source resources to support reproducibility and further research.

%\keywords{
%{OR in sports, Markov decision process, golf strategy}
%} 

\section{Introduction}\label{intro}

Golf is a game in which players aim to put a ball into a \emph{cup} (sometimes also called a hole) with the help of \emph{clubs} (the main types are woods, irons, and putters) using the least number of shots (see \cite{rulesGolf} for the 2019 official rules of golf).  The field where the golfer plays is called a \emph{(golf) course}. It consists of eighteen independent (and different) \emph{holes}. A \emph{hole} comprises different areas as shown in Fig. \ref{golf_hole} (see subsection \ref{sec:golf} for more details).

\begin{figure}[h]
\centering\includegraphics[trim=0 50 0 0 ,width=0.7\textwidth]{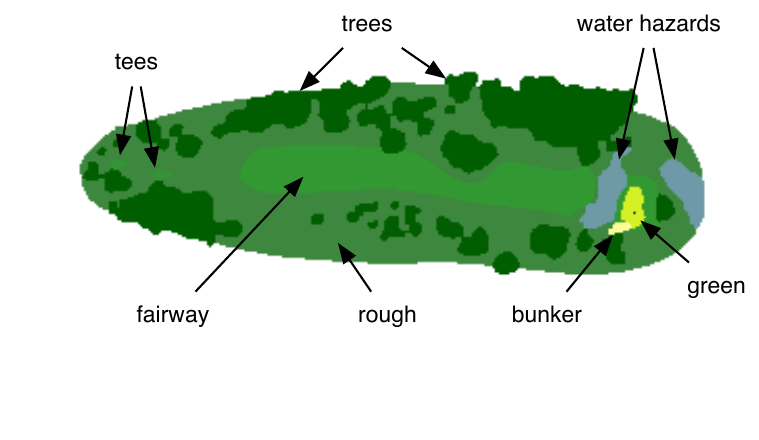}
\caption{A hole and its different areas, the area beyond is out-of-bounds}\label{golf_hole}
\end{figure}

The golfer’s score on a hole is essentially the number of shots taken to put the ball into the cup, plus any possible penalties (when the ball ends up in the water or goes out-of-bounds) \cite{rulesGolf}. The main type of competition is  \emph{stroke play} where the players usually play 4 \emph{rounds} of 18 holes.  The player’s final score is obtained by summing up the score of the 72 holes and the winner is the player with the lowest total score. We focus on this type of competition. Before a competition, elite amateurs and professional players typically practice on the tournament's golf course to identify the main hazards and risks. This usually includes inspecting the course for obstacles, hazards, slopes, and weather conditions. They typically report such information on a golf course booklet and they include personal recommendations on which club to use and target to aim for under different scenarios (pin position, wind, etc.). The main purpose of the course booklet is to help golfers develop an effective strategy for playing their game during the competition. An ideal personalized booklet would describe which shot to play from any (reachable) position on the golf course, given the current performances of the player. Building such a booklet is beyond human capabilities, and in this paper, we develop a methodology to automate the construction of such a  ``strategic'' booklet through the use of (available) historical data on the performance of players and Markov Decision Processes and we show that it is computationally tractable. The underlying problem is strategy optimization.

Markov chains and Markov decision processes are the models of choice for performance assessment and optimization in sports, as they capture the inherent probabilistic nature of the success of every ‘action’ performed by athletes or teams. They have been used, for instance, in tennis, basketball, volleyball, ice-hockey, golf, soccer, darts, and snooker, e.g.~\cite{Terroba2013,trumbelj2012,Routley2015,Pfeiffer2010,Hoffmeister2015,Hoffmeister2017,Heiny2014,Maher2012,BroadieKo,Sugawara}.  Building such a model for golf strategy optimization requires detailed information about the golf course and the player's past performances. In this work we use simple 2D information on golf course extracted from aerial views and historical data on the player's performances taken from the ShotLink database.

The introduction of the ShotLink intelligence program (an initiative of the US Professional Golfer Association (PGA) to share data it collects in real-time on all shots taken on the PGA Tour - the PGA championship - since 2001) has stimulated a lot of academic research  in the past 15 years. Broadie’s pionering work on the stroke-gained method  \cite{Broadie2008,Broadie2012} has revolutionized the analyses of professional golfers' performances on the PGA Tour. In addition, a large body of work exploits the ShotLink database to study various aspects of the game of golf (such as the effect of luck, pressure on performance, the existence of the hot hand phenomenon) through statistical analyses (e.g.~\cite{Baugher2016,Ozbeklik2017,Fearing2011,Connolly2012,Robertson2014,Stockl2018,Gnagy2015,Connolly2008,Connolly2009,Connolly2012b,Hickman2015,Hickman2019,Arkes2016,Heiny2012}), performance prediction through machine learning (e.g.~\cite{Hucaljuk2011,Huang2010,Moorthy2013,Wiseman2016,Lim2017,Drappi2018}, see \cite{Bunker2019} for a recent survey), and the evaluation of different parameters (distance, dispersion, hole size) on performance through simulation and/or optimization \cite{BansalBroadie,BroadieKo}.  This work is part of the area of research called golf analytics. 

Golf strategy optimization was addressed first by \cite{Sugawara}. The authors show how to approximate the optimal  strategy of a player using a skill-model, simulation and Q-learning. The underlying simulation and skill-models are similar to \cite{BroadieKo}, who used a greedy approach to model a golfer's strategy (basically they assume that a golfer always chooses the best shot assuming he will play a perfect one - which is actually a strategy often chosen by amateurs). \cite{Sugawara} show-case their approach using different types of ``average'' players whose skills are characterized by parametrized distributions (the parameters are taken from statistical information available on ``average'' players). In this work, we use similar simulation and skill models as those presented in \cite{Sugawara} and \cite{BroadieKo}, but we apply the methodology using empirical distributions of PGA Tour players, constructed from historical data available through the ShotLink database. 

 We show in particular that the natural Markov Decision Problem associated with the strategic optimization problem (this is essentially the same underlying MDP in \cite{Sugawara}) can be solved exactly in a reasonable amount of time.  In addition, we illustrate how the corresponding methodology can be used for skill improvement. The corresponding approach should help PGA tour professionals (and amateurs collecting similar data using systems such as Arccos) to substantially improve their performances.   

The remainder of the paper is structured as follows. In Section~\ref{method}, we describe how raw ShotLink data are processed to build statistical models of PGA Tour players' skills. Section~\ref{sec:3} presents the simplified physics and rule implementations used to simulate ball trajectories and model the game of golf. Section~\ref{ssp} introduces the stochastic shortest path formulation used for strategy optimization. Section~\ref{sec:6} demonstrates how the framework can be leveraged to evaluate the value of specific golf skills through targeted interventions, and we challenge in particular the ``Drive for show, putt for dough'' adage. Finally, Section~\ref{concl} concludes the paper and outlines directions for future research.

\subsection*{Some golf terminologies} \label{sec:golf}

The \emph{tee} is the area where the player starts (there might be several potential areas but the official starting point for a (round of a) tournament is delimited by \emph{tee balls} - usually corresponding to a rectangular area of around 5\(m^2\)): the grass is short and this is the only place where the golfer can use a tee (same name but referring to a small t-shaped piece of wood) to raise the ball from the ground (and ease the shot). The \emph{green} is the area where the grass is very closely mown, and the cup and \emph{pin} (a flagstick that makes the cup visible from a long distance) are placed. The \emph{fairway} is the part of the hole between the tee and the green where the grass is kept short and shows the (intended) path to the hole. The \emph{rough} is an area of higher grass around the fairway. Usually, the further you get from the fairway, the higher the grass. There are usually different types of \emph{roughs} that vary in grass height and density, namely light or heavy roughs. The \emph{bunkers} are hollows filled with sand that serve as ``traps'' from which it might be difficult to escape (depending on the depth and texture of the sand). The \emph{water hazards} are typically ponds or other bodies of water where the player cannot usually play (unless very shallow or the ball lies on the shore) with a 1-shot penalty to get out (usually close to the point of entrance). \emph{Out-of-bounds} is an area where it is forbidden to play: the player receives a 1-shot penalty if shooting a ball out-of-bounds and has to place the ball in the previous position. There are other ``obstacles'' (formally speaking, obstacles refer to water or bunkers only in the game of golf) such as bushes or trees. See again Fig. \ref{golf_hole} for an illustration of the different areas.

\section{Modelling PGA players' skills from the ShotLink data}\label{method}

In golf, like in many skill games, the result of a shot might differ from the intention. There are several elements that may influence the deviation of a shot from the intended target:   subtle differences in golf swing that can result in different launch parameters (speed, spin, angle, etc.), weather conditions (wind, humidity, etc.) and lie conditions (fairway, rough, bunker, uphill, downhill, buried ball, ball in a divot, etc.). At the strategic level, one needs not take into account all these aspects. Some can be taken into account at the operational level when playing the game. We focus here on the variation of launch parameters and different initial surface (fairway, bunker, rough). To understand the effect of different launch conditions from the same surface (fairway), we present in Figure \ref{trackman} data collected for an elite amateur golfer using TrackMan\footnote{``TrackMan is a radar system that uses Doppler technology to track and record 3D characteristics of a sports ball in motion''. (Wikipedia)}.

\begin{figure}[h!]
\begin{tabular}{cc}
\includegraphics[trim=19cm 0 0 0,clip,width=0.3\textwidth]{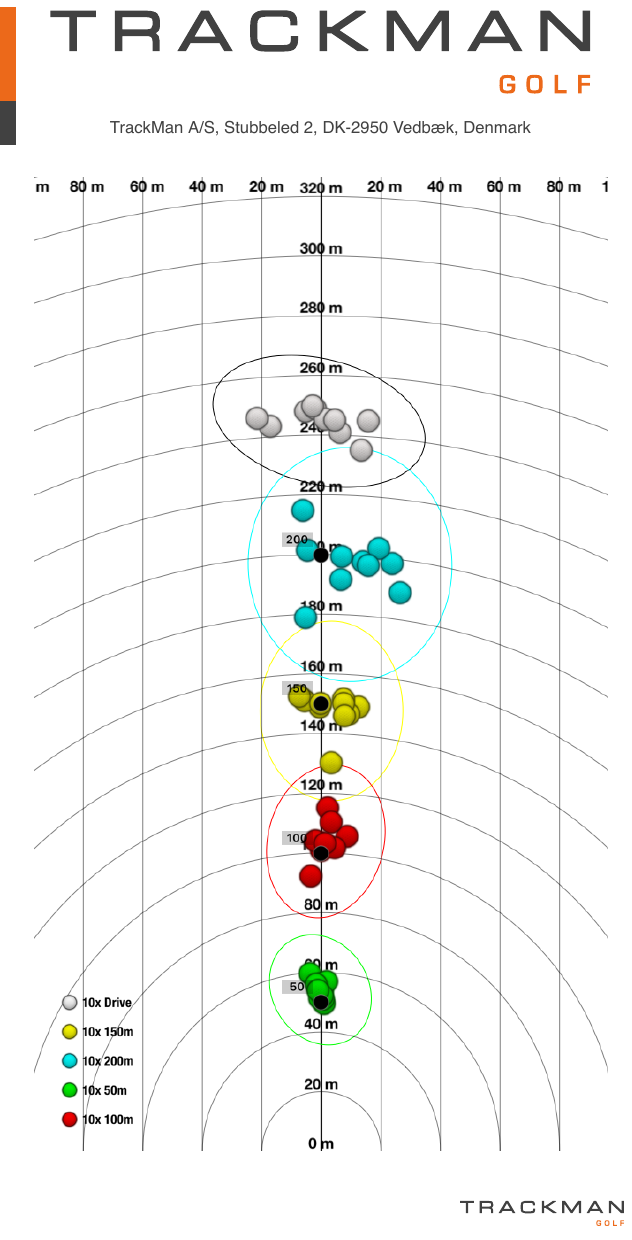} \\
\end{tabular}
\caption{The figure presents a set of target points and the empirical distribution of a player's shots around those targets. The data comes from a member of the French U21 amateur elite team in 2016, with all shots hit from the fairway. We refer to this as a TrackMan profile. The 10 grey points represent drives hit with a driver, demonstrating an average distance of just over 240 meters, along with the distribution of the shots around this average. The blue (for 200 meters), yellow (150 meters), red (100 meters), and green (50 meters) points illustrate the empirical shot distributions around targets located at those distances directly ahead. The ellipses around the points represent a 95\% confidence region, calculated by TrackMan, indicating where most of the shots are expected to land based on the observed data.}\label{trackman}
\end{figure}

A TrackMan profile provides, for certain explicit targets, and given the type of surface, an empirical distribution of the realization of the shots for the player (note that this does not take into account the roll of the ball and this assumes 2D trajectories). These profiles are what we consider in the remainder of the manuscript. More formally, we assume that the position of the balls in 2D, around a targeted point $(0,d)$, with $d\in [0,D_{\mbox{max}}]$, and from a surface $s$ (fairway, bunker, rough) is a random variable $X_{d,s}$ that can be described by a probability density function in 2D  ($D_{\mbox{max}}$ is the maximum distance a player can target). More formally, we consider the sample space $\Omega:=\{(x,y) \in \mathbb{R}^2\}$ and we assume that $X_{d,s}$ follows a probability density function $f_{d,s} : (x,y) \in \Omega \mapsto f_{d,s}(x,y) \in \mathbb{R}_+$ for any $d\in [0,D_{\mbox{max}}]$ and for any $s\in \{fairway, rough, bunker\}$.  The TrackMan profile from Figure \ref{trackman} hence represent samples from the random variables associated with the 5 distances targeted, from the fairway. In this work, we do not impose specific distributional forms (e.g., normal distributions) for the functions $f_{d,s}$. Instead, we rely on empirical distributions derived from observed data, assuming these are representative samples of the corresponding random variables across specific distances and surfaces. Note, however, that parametrized distributions can also be incorporated within our framework through appropriate sampling techniques.

 These empirical distributions can be represented by a TrackMan profile.  Hence in the following, we will consider that skills of the player, outside of the green, are given to us through a TrackMan profile.

Collecting accurate TrackMan profiles of a player is nearly impossible, as it requires the player to hit thousands of balls from the same surface with different target distances. Instead, we will infer approximate TrackMan profiles by exploiting the ShotLink database and common strategies of PGA tour players. The ShotLink database collects the positions of the ball (3D coordinates) of PGA Tour golfers since 2004 in every PGA Tour competition (as well as other information). Recovering TrackMan profiles from such data is not straightforward, as we have no information about the player’s intention: knowing the final destination does not help in assessing the deviation from the intended target. Moreover, we do not have information on the carry and roll of the ball. However, there are some invariants in professional game plans, and we build upon common strategies that professionals use to infer the intention and the ball’s trajectory from the database. Next, we explain the core idea of our approach. It is important to note that our goal is not to replicate exact TrackMan profiles of the players, but rather to create profiles that are sufficiently realistic to effectively demonstrate our methodology on data representative of PGA tour players.

First, professionals tend to target the pin whenever possible and not too risky, which is the case for reasonably short shots - say less than 150m (in reality, they tend to aim slightly off the pin if there is an obvious obstacle close to the pin (e.g.~a water hazard or bunker), or if the green is not flat, as they generally prefer uphill to downhill putts, but we omit such details). For longer shots, they might play it a bit safer, aiming between the pin and the middle of the green, but since we do not have information about the geometry of the green (and the position of the pin with respect to this geometry), we also assume in this case that they target the pin again.  Second, professional golfers tend to choose a general strategy for each of the tee shots (on par-4 or par-5) before the first day of the tournament (in the training rounds on the previous days). So unless there are very different conditions between two rounds (e.g.~tee positions, weather conditions), we can reasonably assume that they aim at the same target.

As mentioned, we do not have information about the carry and roll of the ball (nor on the potential lateral spin). We therefore assume for simplicity - and lack of data - that the trajectories are straight and that the final endpoint is independent of the carry/roll trajectory, so that the empirical distribution of a shot around the target only depends on the distance and on the lie of the ball. We believe that the first assumption is not very restrictive, as most shots played by professional golfers are straight (some players may have a slight preference for a certain lateral spin - fade or draw - but this is usually very subtle). For the second hypothesis, the main shortcoming might derive from very short game situations - say below 30m - where different trajectories can be chosen (with a preference usually for rolling the ball on the green as much as possible, but where a ``lob shot'' is needed when there is an ``obstacle'' (say a bunker) close to the pin and in between the ball and the hole). In such situations, we implicitly assume that the distribution of the shots around the target does not depend on the type of shot played, which is certainly a limitation. For other shots, the ball tends to roll more on the fairway and the green than in the rough when coming from a long distance (but not much when at a distance of less than 150m). However, we believe that this effect is limited.

So, in concrete terms, we build the TrackMan profiles from the trace of the shots in ShotLink. The position of the ball is stored as a 3D coordinate $(x,y,z)$ in the database, where $x$ and $y$ are essentially the longitude and latitude of the ball (possibly expressed in a different coordinate system) and $z$ is the elevation. We assume here for the exposition that the coordinates are expressed in meters. We get rid of the $z$ coordinate as we assume, in this project, that the trajectories are flat, for computational efficiency. We proceed as follows. Assume that $(x_0,y_0)$ is the position of the ball before the shot and that the ball lies on the fairway, $(x_1,y_1)$ is the position of the ball after the shot, and $(x^P,y^P)$ is the position of the pin.  Let $d=||(x-x^P,y-y^P)||_2=\sqrt{(x-x^P)^2+(y-y^P)^2)}$ and let $M$ be the (transpose of the) rotation matrix, that maps $(x^P-x_0,y^P-y_0)$ to $(0,d)$. We assume that $M \cdot (x_1-x_0,y_1-y_0)$ is drawn from $X_{d,s}$ (with $s=fairway$) and thus we can build proxies of TrackMan profile using this procedure by using all data from a player in the database (we restrict to data from the last 12 months). For tee shots, we assume, as explained above, that the player targets the average position over the 4 rounds that he played (if he played only two rounds, we discard the corresponding data - we also check that the tee positions are within a 5 meter radius over the 4 rounds) and we applied a similar procedure.

We have inferred the TrackMan profiles from this procedure. The Fig. \ref{fig:philfairway} ``original'' panel shows the outcome for Phil Mickelson on fairway shots. Two things clearly emerge from the Fig. \ref{fig:philfairway} ``original'' panel (or Fig. \ref{fig:philrough} and Fig. \ref{fig:philbunker} for shots played from the bunker and the rough, see Appendix). First, we do not have historical data for any distance. Second, there are pairs that look suspicious: the pair with a target point at a distance of roughly 370m and final destination at 250m is probably an outlier. Indeed in this case, it is fairly obvious that the player did not try to target the pin (way too far to be within reach), but probably played a safe shot on the fairway between his position and the target. One way of detecting this kind of outlier is to put a cap on the maximum distance error that a player could make in principle. Indeed, professional players are known to be fairly accurate in terms of distance control (unless there are very particular conditions, e.g.~hard ground, wind blowing suddenly, hit of a tree, etc\ldots). We used the TrackMan data for PGA Tour professionals taken from https://blog.trackmangolf.com/category/tour-stats (see Fig. \ref{fig:philfairway} (left)) to set the cap, as we explain below.

We chose to keep all wedge shots (that is, shots under 100m), as in this case, the assumption of targeting the pin makes perfect sense (there might still be outliers, for instance, if the player’s ball ends up in a water hazard and the player drops it close to the entry point, then the observed deviation does not correspond to the true one), while we cap the maximum distance error to 20m  for shots over 100m (the maximum in the PGA Tour TrackMan data available to us is 15m, see Fig. \ref{fig:philfairway} (left)). The results are presented in the ``cleaned'' panel of Fig. \ref{fig:philfairway} . This choice is consistent with statistics from the literature (see Fig. 3 in \cite{James2008}), and we validated this number with two trainers of the French U21 elite amateur team using data from their elite amateurs. Fig. \ref{fig:someplayer1} shows the results for 24 PGA Tour professionals.

\begin{figure*}[h!]
\includegraphics[width=0.7\linewidth]{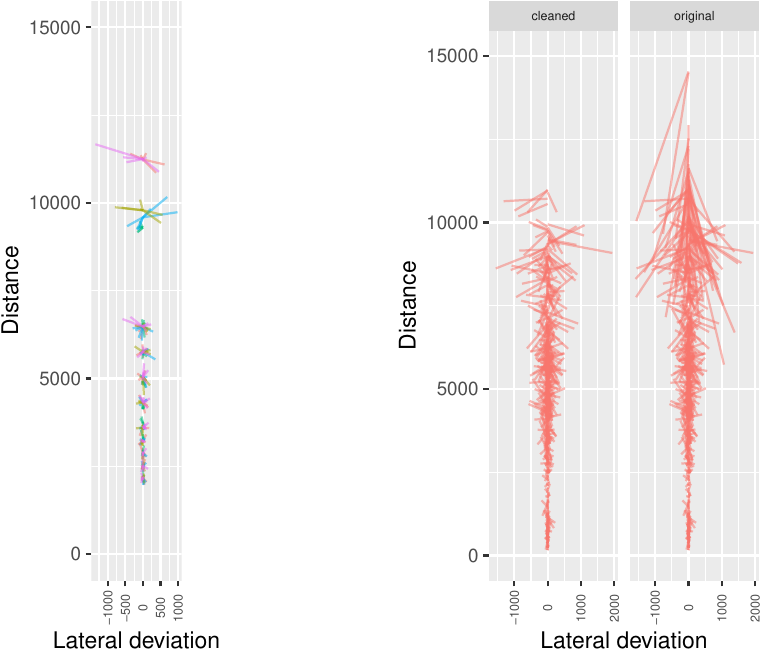} \caption{\label{fig:philfairway} Each segment of the figure represents a target/destination pair for shots played from the fairway. All target points have been rotated so as to appear on the y-axis. The ``original'' panel shows inferred data for Phil Mickelson before removing outliers. The left figure shows the (few) TrackMan data we have for PGA Tour professionals (the different colors represent different players). The ``cleaned'' panel shows the data for Phil Mickelson after the removal of outliers. Numbers are in meters.}\label{fig:unnamed-chunk-1}
\end{figure*}

\begin{figure*}[h!]
\includegraphics{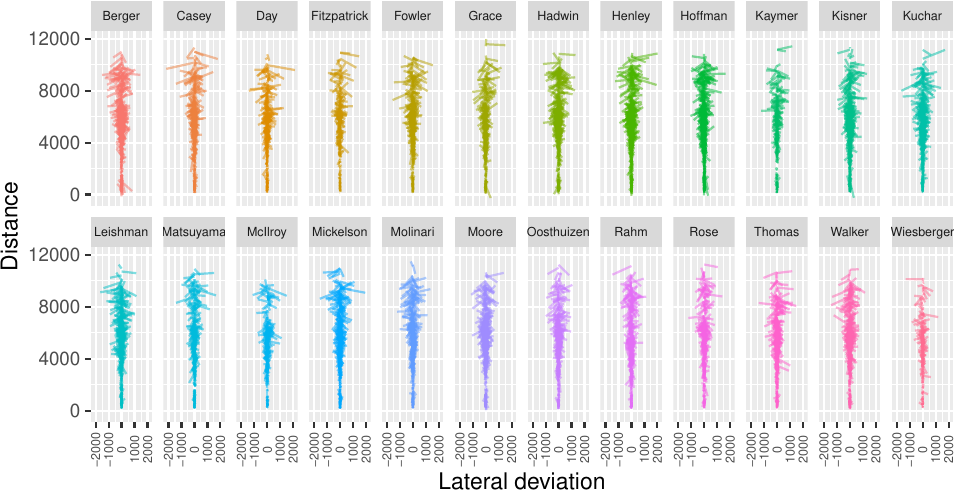} \caption{\label{fig:someplayer1} Each segment of the figure represents a target/destination pair for shots played from the fairway after removal of outliers. All target points have been rotated so as to appear on the y-axis.  Numbers are in meters.}\label{fig:unnamed-chunk-2}
\end{figure*}

Although rigorous statistical validation is not possible (since we do not have the true TrackMan profiles for the corresponding players), the figures are fairly consistent with the limited TrackMan data we have (see Fig. \ref{fig:philfairway} (left)). However, there are still some obvious outliers (e.g.~Grace and Kisner have shots that end up behind them - perhaps due to hitting a tree), but we are not too concerned about these few remaining cases, as the effect of these outliers will be smoothened in the next preprocessing phase.

The task is somewhat harder from the rough, as there are many types from light to heavy, and the lie of the ball might vary substantially and might have a strong impact on the outcome: it could be more difficult to hit a buried ball in a light rough than a ball lying on the surface of a heavy rough. As we do not have access to such information, we collected data from the rough without differentiating these situations. Of course, this is a limitation. For similar reasons as the fairway shots, we can remove outliers when the distance error is too large. We do not have TrackMan data of PGA professionals from the rough, so we need to set the thresholds somewhat more arbitrarily. The distance control error might be more important in the rough. The accuracy error on fairway and rough shots was investigated in \cite{James2008}. Based on the statistics reported in Fig. 2 of the corresponding manuscript, we set the threshold on distance control error for outlier detection at  30m. We again validated this number with two of the trainers of the French U21 elite amateur team using data from their elite amateurs. The results are shown in Fig. \ref{fig:philrough} and Fig. \ref{fig:someplayer1rough} (see Appendix).

For the bunkers, we apply the exact same strategy as for the rough. The results are shown in Fig. \ref{fig:philbunker} and Fig. \ref{fig:someplayer1bunker} (see Appendix). We have very scarce data for distance over 40-50 meters, which is not surprising since there are many more bunkers around the greens than so-called \emph{fairway bunkers}. We will consider how to interpolate the missing data. While interpolation has limitations (especially where we have very few points), it has minimal consequences, as few shots are played from fairway bunkers.

For driving data off the tee, the results are shown in Fig. \ref{fig:philtee} and Fig. \ref{fig:someplayer1tee} (see Appendix). Here again we applied a threshold of 30m  to the distance control error for outlier detection.

To create complete TrackMan profiles, we would need empirical distributions for any possible distance. We use a form of bootstrapping to generate missing data. Observe that this procedure would be needed even if we had true TrackMan profiles from the beginning, such as the one in Fig. \ref{trackman}, as there contain only a limited number of targeted distances and sample points. The main idea is to use a local linear approximation. For any distance \(d\) we might target, the basic idea is to grow a disk around \((0,d)\) until it contains enough sample target-destination pairs (from the inferred data), and scale the coordinates of the arrival points by the ratio of the original targeted distance and \(d\) (that is, assuming a linear relation: if the targeted point was \((0,t)\) (within the disk) and the arrival point \((x,y)\), we assume that for the hypothetical target \((0,d)\) this would have resulted in the arrival point \(\frac{d}{t}\cdot(x,y)\)). The radius of the disk is defined so as to grab enough data to be statistically relevant, but not too many when this includes points that are too remote - this is the case when we have fewer data available, such as for fairway bunkers, for instance: we set the number of points to 50 and put a cap of 30 meters on the radius unless we have fewer than 10 points available (in which case, we take the closest 10 points). We have also assumed that the distribution is symmetric along the y-axis, and that the lateral error and the distance control error are independent. We used these hypotheses to ``shuffle'' the x and y coordinates of the arrival points and avoid too much bias toward the existing data points. Additionally, we used the 95th percentile of the maximum observed target distance as the maximum targetable distance of a surface, and capped the maximum distance that the ball might reach by the maximum observed distance. The results are shown in Fig. \ref{fig:bootstrapfairway1}, and Fig. \ref{fig:bootstraprough1}, Fig. \ref{fig:bootstrapbunker1}, Fig. \ref{fig:bootstraptee1} and Fig. \ref{fig:bootstrapfairwaytee1} (see Appendix). In all the figures, we have restricted the target distances to multiples of 2.5m, and generated 15 realizations for each distance.  Note that for consistency, we ensured that the average lateral dispersion would increase with distance. Hence, we slightly rescaled the lateral dispersion by the inverse of the ratio with the average of the previous distance when this was not the case. Furthermore, we paid attention to the fact that the average lateral deviation from the rough and from the bunker (for a given distance) could not be less than that from the fairway to compensate for missing data that could bias the results. Although there are still a few inconsistent points, the result appear fairly clean. While the parameters used above could be fine-tuned, we are satisfied with the results from the current settings (see Table \ref{tab1} and Table \ref{tab2} and the associated discussion) as, again, our objective is to construct realistic PGA Tour player profiles, not to create exact virtual clones. 

\begin{figure*}[h!]
\includegraphics{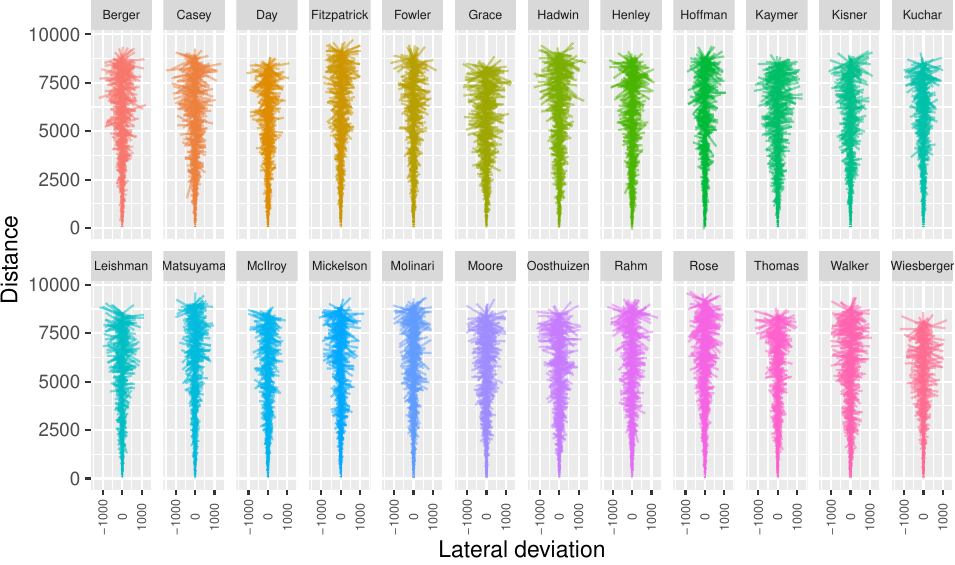} \caption{\label{fig:bootstrapfairway1} Each segment of the figure represents a target/destination pair for shots played from the fairway after bootstrapping. All target points have been rotated so as to appear on the y-axis.  Numbers are in meters.}\label{fig:unnamed-chunk-9}
\end{figure*}

Now the last skills we did not consider yet are putting skills. On a green, a professional PGA Tour golfer typically makes 1, 2, or 3 putts (perhaps 4 or 5 putts in exceptional situations, but this represents less than 3 cases out of 10000 in our data, so we assume that only these three situations can occur). As we do not have details about the slopes of the greens, we assume that the average number of putts is simply a function of the distance to the hole, that is it follows a certain function $p:d\in \mathbb{R} \mapsto p(d) \in \mathbb{R}$. In order to evaluate this function, we estimated the probability of 1-putt, 2-putts, and 3-putts for any possible distance on a green. We focus now on the 1-putt and 3-putts probabilities, as the 2-putts probability is easily computed from the other two.

\begin{figure*}[h!]
\includegraphics[width=0.6\linewidth]{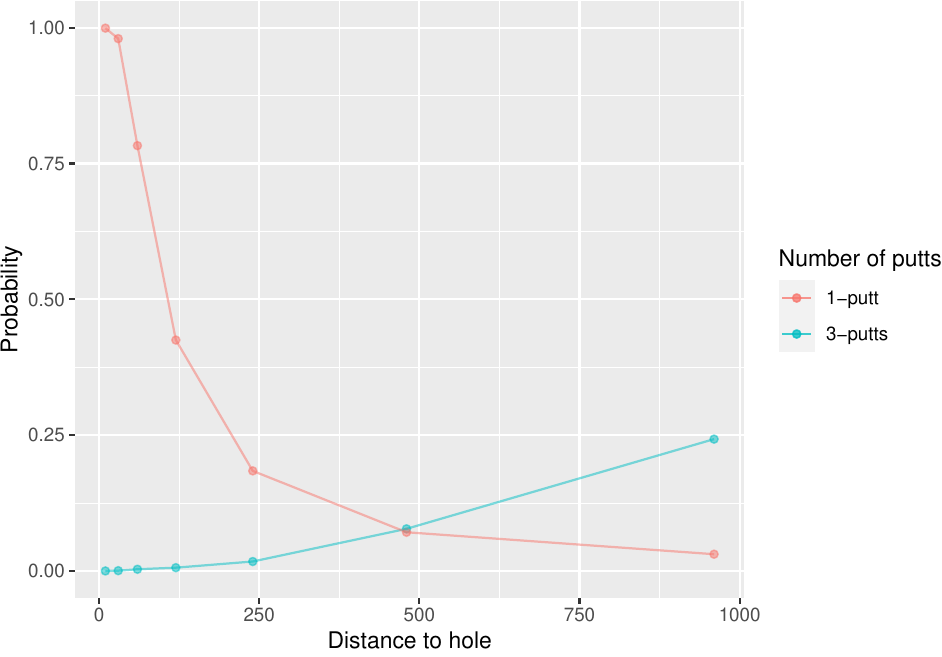} \caption{\label{fig:tourputting} 1-putt (in red) and 3-putt (in blue) probabilities as a function of distance for the ``average'' PGA Tour player. Numbers are in meters.}\label{fig:unnamed-chunk-15}
\end{figure*}

The longest possible distance on a green (that is, the diameter of the geometrical object) is usually no more than 25m, especially in the US (the old course in Saint-Andrews (UK) is an exception with a diameter of some greens reaching 50m). We thus consider distance to hole below 32m (above this value, we have very little data: about 1 per 16,000). To build the histogram, we need enough data for each distance. Unfortunately, we have more data close to the pin than far from it (as players usually get closer each time they play without putting the ball into the hole). For individual players, we thus inferred the probability from the data available below 16m. We used buckets of doubling size (to have enough data within each bucket for statistical relevance - more than 30 points typically) with the following ``breakpoints'' (in meters): (0, 0.5, 1, 2, 4, 8, 16) and assigned the value to the midpoint of the interval. so, in concrete terms, for a given player, we collected all data regarding putts made from a distance in the range of distance, say $[1,2]$, and we estimated the 1-putt (2-putt, 3-putt resp.) probability as the frequency of 1-putt (2-putt, 3-putt resp.). We then assumed that all putts were played from the mid-distance, that is 1.5 m here. For distances between 16m and 32m, due to the lack of data for individual players, we aggregated all data from all players to estimated the corresponding probabilities. We then used a linear interpolation to build a proxy of the function $p$ for any player.

The results are presented in Fig. \ref{fig:puttingexamplesproba}. If we aggregate all the data from professional PGA Tour players, we obtain the result in Fig. \ref{fig:tourputting}, which is in line with the literature (see Fig. 1 in \cite{James2008}). We also compare the average number of putts as a function of the distance for different PGA Tour players. The results are shown in Fig. \ref{fig:puttingexamplesaverage} (of course, our estimations could probably be improved and smoothened, but we believe that using more advanced approaches is not relevant at this stage, as there are other simplifications in our models that probably dominate this one).

We preprocessed the data with R, and the corresponding code is available, upon request, in a companion zip file for the replicability of our results.

\begin{figure*}[h!]
\includegraphics[width=0.6\linewidth]{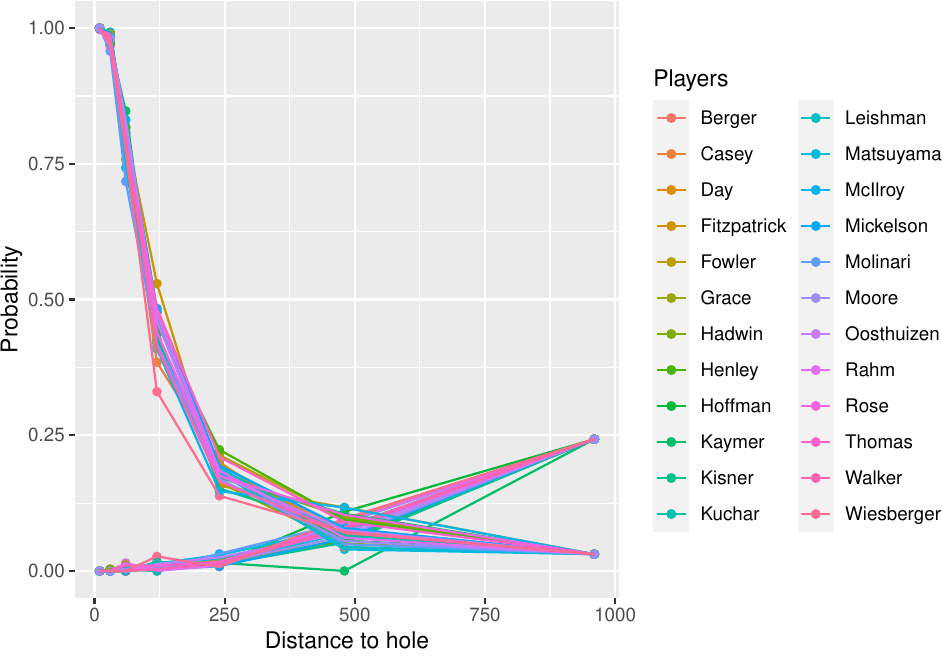} \caption{\label{fig:puttingexamplesproba} 1-putt and 3-putt probabilities as a function of distance for a subset of players. The 1-putt probabilities typically appear as non-increasing curves, while the 3-putt probabilities generally follow non-decreasing trends. Numbers are in meters.}\label{fig:unnamed-chunk-16}
\end{figure*}

\begin{figure*}[h!]
\includegraphics[width=0.5\linewidth]{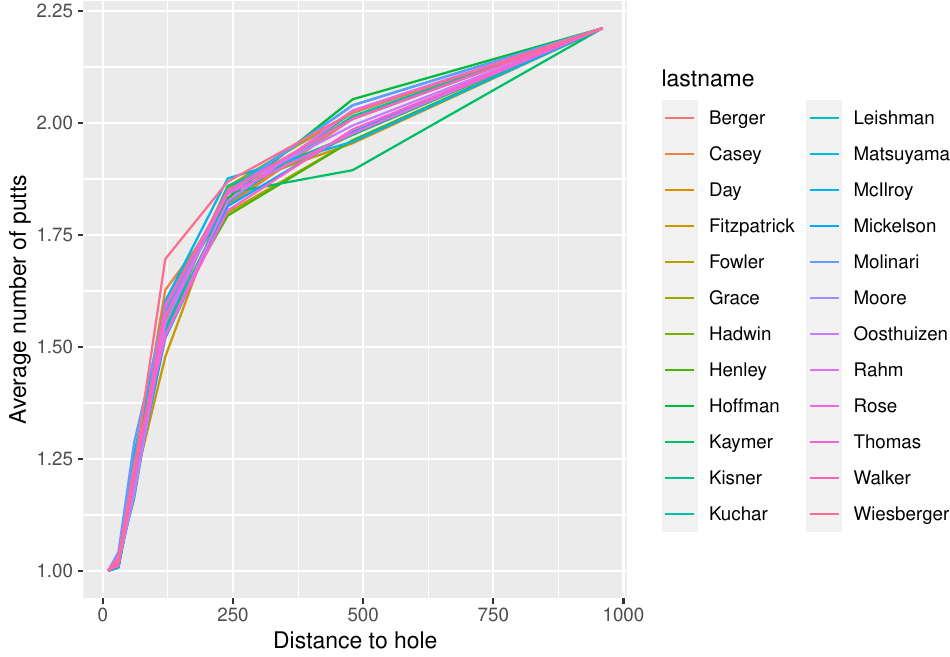} \caption{\label{fig:puttingexamplesaverage} Average number of putts as a function of distance for a subset of players. Numbers are in meters.}\label{fig:unnamed-chunk-17}
\end{figure*}

\section{Modelling the game of golf}\label{sec:3}

In order to optimize a golfer strategy with a Markov Decision Process, we need a model to predict the (stochastic) outcome of any ``action'' we may use. In our case, as will become clear later, actions correspond to shots and we need to specify the result of any (possible) shot. 

We explained in the previous section how we inferred reasonable TrackMan profiles of players on each surface and for each distance given the ``trace'' of the players on the different tournaments available in the database. These profiles are 2D and we assume in the following that they represent the projections of the \emph{flight} of the ball: we thus implicitly assume that the ball does not roll and that the trajectories are straight. Now we use a 2D representation of the golf course as well, using stylized 2D pictures of the holes and a clear encoding of each surface and obstacles similar to Fig. \ref{golf_hole} (bushes/trees in dark green, roughs in green, fairways in light green, greens in yellow green, and bunkers in egg nog). We actually created (manually) 2D raster of the holes of three golf courses using aerial views from google maps: Augusta National Golf Club (that hosts one of the four major tournaments: the Masters), Le Golf National (that hosted the Ryder Cup competition in 2018, a very famous biennial competition between male teams from Europe and the United States) and the Bay Hill Club and Lodge, Orlando (that hosts the Arnold Palmer Invitational). We chose the resolution so that 1 cell roughly represents a square region of side 1m (typically between 0.7m and 1.5m, depending on the hole).

For a position $(x,y)$ on a surface $s$ (actually a cell in the raster), a shot is essentially a selection of a target point $(x^t,y^t)$. This target point can be characterized by a distance $d$ and an direction/angle.  Let $M$ be the (transposed of the) rotation matrix associated with the corresponding angle. The target point is $(x,y) + M\cdot (0,d)$. Now in order to estimate the distribution of the outcome of the shot, we consider $k$ realizations sampled from $X_{d,s}$. Let $(x_1,y_1) ,...,(x_k,y_k)$ be the corresponding samples. In the absence of trees, water hazard and out-of-bound area, the distribution of the outcome could be approximated by the sample $(x'_1,y'_1) ,...,(x'_k,y'_k)$, where $(x'_i,y'_i)= (x,y) + M\cdot (x_i,y_i)$ for all $i=1,...,k$. We call this the hypothetical empirical distribution. 

In the presence of trees (and the like) in the ball’s trajectory, we assume that the ball will stop on the trajectory right before the first obstacle it encounters. That is, the ball will hit the obstacle and it will neither ``bounce off'' nor penetrate the obstacle, it will simply stop (consider a collision with a dense tree, like fir).  More precisely, we assume that trees are infinitely high, and that the ball simply stops right before the contact point.  When the ball falls into a water hazard, we assume that the player will ``drop'' the ball at the entry point (which is the most common option out of the different possible options for a golfer in such situation). When the ball ends up out-of-bounds, we (re)position the ball at the origin of the shot (with a 1-shot penalty according to the rules of golf). 

Technically speaking, we apply the following procedure to identify the destination cell. We use Bresenham's algorithm  (\cite{Bresenham}) to identify an ordered set of cells from the raster/picture that are traversed by the trajectory. Suppose we consider the trajectory from $(x_1,y_1)$ to $(x'_1,y'_1)$. Let $\mathcal{C}=\{c_1,...,c_l\}$ be the ordered set of cells return by Bresenham's algorithm. Let $i$ be the index of the first tree cell in $\mathcal{C}$ ($i=l+1$ if there is none). We first truncate the trajectory to $\mathcal{C}'=\{c_1,...,c_{i-1}\}$ (observe that we do not allow a player to play from a tree so $i>1$). Then, if $c_{i-1}$ is a water cell, we let $j$ be the largest index such that $c_j$ is not a water cell (again here we do not allow a golfer to play from a water hazard so $i-1>j>1$) and we truncate the trajectory to $\mathcal{C}''=\{c_1,...,c_{j-1}\}$. Finally if $c_{l}$ is an out-of-bound cell, we truncate the trajectory to $\mathcal{C}''=\{c_1\}$. An illustrative (non-realistic) example of a simulation is given in Fig. \ref{figsimulation}.

\begin{figure}
\begin{tabular}{ccccc}
\includegraphics[width=0.2\textwidth]{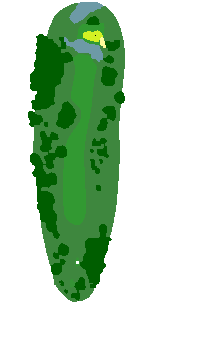}&
\includegraphics[width=0.2\textwidth]{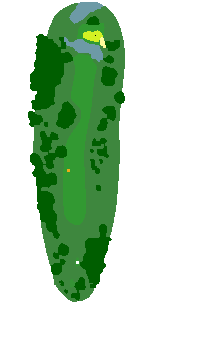}&
\includegraphics[width=0.2\textwidth]{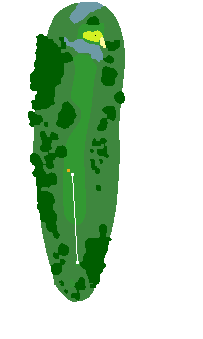}&
\includegraphics[width=0.2\textwidth]{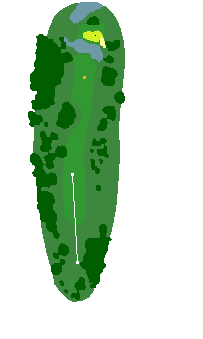}&
\includegraphics[width=0.2\textwidth]{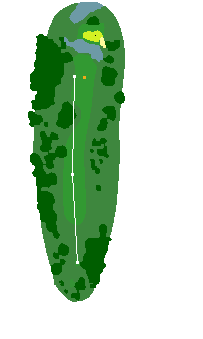}\\
\includegraphics[width=0.2\textwidth]{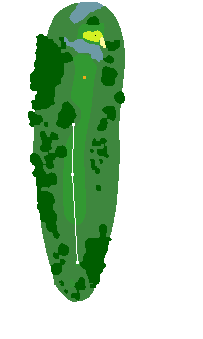}&
\includegraphics[width=0.2\textwidth]{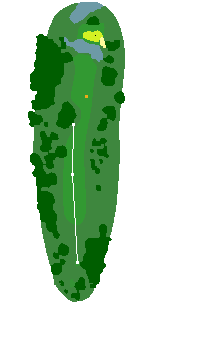}&
\includegraphics[width=0.2\textwidth]{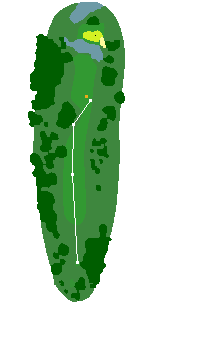}&
\includegraphics[width=0.2\textwidth]{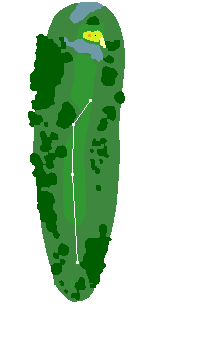}&
\includegraphics[width=0.2\textwidth]{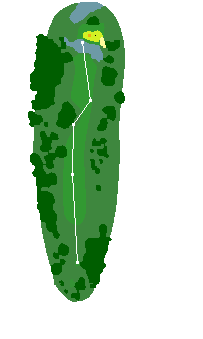}\\
\includegraphics[width=0.2\textwidth]{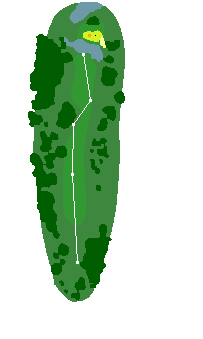}&
\includegraphics[width=0.2\textwidth]{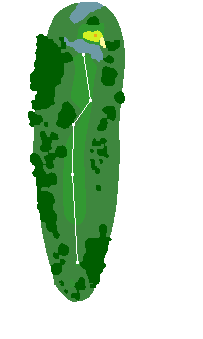}&
\includegraphics[width=0.2\textwidth]{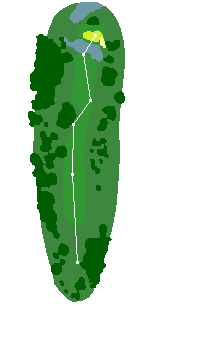}\\
\end{tabular}
\caption{A sequence of simulated shots (pictures numbered sequentially from left to right and top to bottom): the white dot in picture 1 corresponds to the initial tee position, the orange dot corresponds to the intended target, white segments show the player’s trail (note that the first time a segment appears, it reflects the realization according to the hypothetical empirical distribution - without the detection of collisions and special events: it might be shortened in a second stage, e.g. pictures 6 and 11, taking into account the obstacles, water hazards, and out-of-bounds). Note that we stop as soon as we reach the green, since we can then generate the number of putts as described above. In such situation, assuming that the player makes 1 putt on the green, he would score 7 on the hole (5 shots + 1 penalty - as one shot ended up in the water hazard - to reach the green and then 1 putt). This data is entirely artificial and unrealistic for a PGA professional, and is provided solely for illustrative purposes.}\label{figsimulation}
\end{figure}

We have implemented the corresponding model in C++ for better performance (as will become clear later, we need to call this model a billion times just to create the stochastic shortest path model). The corresponding C++ code is available, upon request, in a companion zip file for the replicability of our results.

\section{The optimization model} \label{ssp}

Before explaining the model, we start with a brief introduction of the stochastic shortest path problem, following  \cite{GUILLOT2018}.

The stochastic shortest path (SSP) problem is a Markov decision process (MDP) that generalizes the classic deterministic shortest path problem. We want to control an agent who evolves dynamically in a system composed of different \emph{states}, so as to converge to a predefined \emph{target}. The agent is controlled by taking \emph{actions} in each time period (we focus here on discrete time (infinite) horizon problems): actions are associated with costs, and transitions in the system are governed by probability distributions that depend exclusively on the previous action taken, and are thus independent of the past. We restrict to finite state/action spaces: the goal is to choose an action for each state, i.e. a \emph{deterministic and stationary policy}, that reaches the target state with probability one (such a policy is called proper), so as to minimize the total expected cost incurred by the agent before reaching the (absorbing) target state when starting from a given initial state. The problem is well-defined when there is a way to reach the target from any state and when there is no improper policy that allows accumulating an infinitely negative cost \cite{GUILLOT2018}.

More formally, a \emph{stochastic shortest path} instance is defined by a tuple \((\mathcal{S}, \mathcal{A},J,P,c)\) where \({\mathcal S}=\{0,1,\ldots,n\}\) is a finite set of \emph{states}, \({\mathcal A}=\{0,1,\ldots,m\}\) is a finite set of \emph{actions}, \(J\) is a 0/1 matrix with \(m\) rows and \(n\) columns and general term \(J(a,s)\), for all \(a\in \{1,...,m\}\) and \(s\in \{1,...,n\}\), with \(J(a,s)=1\) if and only if action \(a\) is available in state \(s\), \(P\) is a \emph{row substochastic matrix} (a \emph{row substochastic} matrix is a matrix with nonnegative entries so that every row adds up to at most $1$. Observe that it is not a usual stochastic matrix as state $0$ and action $0$ are left out) with \(m\) rows and \(n\) columns and general term \(P(a,s):=p(s|a)\) (probability of ending in \(s\) when taking action \(a\)), for all \(a\in \{1,...,m\}\), \(s\in \{1,...,n\}\), and a cost vector \(c \in {\mathbb R}^{m}\). The state \(0\) is called the \emph{target} state and the action \(0\) is the unique action available in that state. Action \(0\) leads to state \(0\) with probability \(1\). We denote with \({\mathcal A}(s)\) the set of actions available from \(s\in \{1,...,n\}\) and assume without loss of generality ({ If not, we simply duplicate the actions}) that for all \(a \in \mathcal{A}\), there exists a unique \(s\), such that \(a \in \mathcal{A}(s)\). We denote with \({\mathcal A}^{-1}(s)\) the set of actions that lead to \(s\), i.e.~\({\mathcal A}^{-1}(s):=\{a: P(a,s)>0\}\).

A (deterministic and stationary) \emph{policy} \(\Pi\) is a function \(\Pi: s \in {\mathcal S}\mapsto {\mathcal A}(s)\),
that is, it assigns an action for each possible state. Let \(y^\Pi_k \in {\mathbb R}_+^n\) be the substochastic (in general, not a purely stochastic vector, as state $0$ is left out.) vector representing the state of the system in period \(k\) when following policy \(\Pi\) (from an initial distribution \(y^\Pi_0\)). That is, \(y^\Pi_k(s)\) is the probability of being in state \(s\), for all \(s=1,...,n\) at time \(k\) following policy \(\Pi\). Similarly, we denote with \(x_k^{\Pi} \in \mathbb{R}_+^m\) the substochastic (in general, not a purely stochastic vector, as action $0$ is left out) vector representing the probability of performing action \(a\), for all \(a=1,...,m\), at time \(k\) following policy \(\Pi\). Given a policy \(\Pi\) and an initial distribution \(y^\Pi_0\) at time \(0\), by the law of total probability (and because each action is available in exactly one state), we have \(x_k^{\Pi}=\Pi^T \cdot y^\Pi_{k}\) for all \(k\geq 0\).

Given a state \(s\in \{1,...,n\}\), a policy \(\Pi\) is said to be \emph{$s$-proper} if \(\sum_{k\geq 0} x_{k}^{\Pi}\) is finite, when \(y^\Pi_0:= e_{s}\). Observe that \(\sum_{k\geq 0} y_{k}^{\Pi}\) is also finite for s-proper policies (as \(y_{k}^{\Pi}=P^T x_{k-1}^{\Pi}\)). In particular, \(\lim_{k\rightarrow +\infty} y^\Pi_{k}=0\), and thus the policy leads to the target state \(0\) with probability \(1\) from state \(s\). An \(s\)-proper policy is thus a policy that converges to the target with probability one, and whose expected number of visits in each action is finite. The expected cost of such policy is thus the well-defined value \(c^T \sum_{k\geq 0} x_{k}^{\Pi}\). The \emph{$s$-stochastic-shortest-path problem }(\(s\)-SSP for short) is the problem of finding an \(s\)-proper policy \(\Pi\) of minimal cost \(c^T \sum_{k\geq 0} x_{k}^{\Pi}\).

As explained in the previous sections, we built (discrete) 2D models for holes and TrackMan profiles, and a 2D simulator of ball trajectories (using Bresenham algorithm). With these elements, we can evaluate a player’s performance for any strategy (a strategy is the choice of shot for any position on the corresponding hole, that is, essentially a choice of direction and targeted distance, according to our 2D representations) on any hole. Indeed, given a strategy, we can build a Markov chain whose states are the possible positions on the hole (the pixels basically), and the transition matrix can be built from the choice of direction and targeted distance in any state as follows: take the different empirical realizations from the TrackMan profile corresponding to the surface where the ball lies, and simulate the outcome of the different realizations on the corresponding hole (we have 15 realizations for each shot from the TrackMan profiles generated in Section \ref{method}). This provides an empirical distribution over the state space and an expected cost for the corresponding ``action'' (1 if no penalty occurs). If the strategy is sound (that is, converging to the target from any position), the corresponding Markov chain is absorbing, and the expected number of steps before being absorbed by the cup can be easily evaluated through computing the fundamental matrix (see for instance \cite{Proba}).

We can go one step further and find the optimal strategy by building an absorbing Markov decision process (an SSP in fact) by simply adding all the possible sets of actions available in a given state to the Markov chain described above. That is, in this SSP model, the states are still the positions on the hole (again the pixels), the actions are the triplets (state, targeted distance, direction), the (empirical) transition matrix and the costs are computed as explained above for the Markov chain model. { We have restricted the set of possible directions to an angle (in radian) in $\{ 0,\frac{2\pi}{180}, \ldots , 179 \cdot \frac{2\pi}{180} \}$ (with this discretization, the player has an aiming precision of at least $\approx$ 1.75m at a distance of 100m), and the targeted distance from the TrackMan profiles is restricted to multiples of 2.5m as discussed earlier).}

In our 2D representations of the hole, we ensured that only locations where the target could be reached were kept (for instance, no rough surrounded by trees or out-of-bounds only: in such case, we would redefine the corresponding zone as a tree area or as an out-of-bounds area), and hence our models are well-defined instances of SSP.

The corresponding instances have the order of 10 thousand states and 150 million actions. We implemented the value iteration algorithm in C++ (see \cite{GUILLOT2018} for details of the algorithm). Most of the time spent on the SSP problem resolution actually entails creating the model (with more than a billion calls to the simulator and Bresenham's algorithm, as we have 15 realizations to simulate for each 150 million actions). Although we have attempted to optimize our code as much as possible, optimizing the computational performance is not the main purpose of this study. Indeed, the computational performance would clearly improve by parallelizing the construction of the model. Instead again we aim at showing that the models are tractable computationally to stimulate further investigations. 

An example of output of the model is provided in Figure \ref{fig:application} (see Appendix).

\subsection*{Computational experiments}

We conducted our experiments on the SD530 nodes within the Curta platform, which consists of 336 nodes. Each node is powered by an Intel Xeon Gold SKL-6130 processor operating at 2.1 GHz, has 32 cores, and is equipped with 96 GB of RAM. More details on the hardware can be found at \url{https://redmine.mcia.fr/projects/cluster-curta}. 

Our analysis focused on 107 golf players from the 2018 Arnold Palmer Invitational (out of 165 participants) for whom we had ShotLink data available to construct TrackMan profiles. We utilized all relevant data from 2017, as well as data from early 2018 leading up to the tournament in March. The average time required to build the model and optimize the strategy for a single hole was approximately 27 minutes, with a standard deviation of 4.5 minutes. The fastest and slowest times recorded were 12 minutes and 44 minutes, respectively. Note that the main reason for setting the number of realizations to 15 was memory limits: the minimum and maximum memory consumption was 67 and 77 GB respectively (and no more than 92 GB could be reserved on the nodes).

The bulk of the computational effort is devoted to model building (a couple of minutes is needed for value iteration typically - 115 seconds on average with a standard deviation of 42 seconds and a maximum of 304 seconds). The  computation time for model creation is essentially multi-linear in the parameters that influence the number of simulations/Bresenham calculations. These parameters include the angle discretization, which affects the number \(d\) of possible shot directions; the target distance discretization, which affects the number \(t\) of target distances; and the number \(r\) of realizations used to generate the bootstrapped TrackMan profiles. Therefore, the empirical computation time is of the order  \(O(d \cdot t \cdot r)\).

Significant improvements in computation time could be achieved through parallel processing. Since the outcomes of individual actions within the model are independent, parallelizing the model creation over \(M\) machines could reduce the time to \(O\left(\frac{d \cdot t \cdot r}{M}\right)\). The value iteration algorithm could also be parallelized, although this would require careful organization due to the dependencies between computations. However, the focus of this paper is to demonstrate the feasibility of solving exactly the SSP model associated with golf strategy optimization, rather than exploring computational optimization through parallelization. We should note also that one could easily reduce the action space by filtering out actions that are dominated, e.g. aiming in a direction in the opposite side of the pin is most of the time clearly suboptimal. 

The corresponding C++ code and script to launch the code on the Curta platform (or similar) are available, upon request, in a companion zip file.

\section{Revisiting the ``Drive for show, putt for dough'' Adage}\label{sec:6}

Traditional performance statistics—such as strokes gained~\cite{Broadie2008,Broadie2012}—help players benchmark themselves against peers but offer limited guidance on {how} to improve. Should a player prioritize distance control, reduce lateral dispersion, work on putting consistency, or increase driving distance? Our simulation framework offers a unique advantage: it enables controlled interventions by selectively substituting individual skill components to estimate their marginal impact on scoring.

To illustrate this capability, we revisit the well-known adage, ``Drive for show, putt for dough.'' Using our simulation environment, we quantify the marginal value of driving and putting by replacing each of the 107 players’ original skill profiles with those of two exceptional performers: Rory McIlroy for driving and Tiger Woods for putting. Figures~\ref{gap_driving} and~\ref{gap_putting} show the resulting distributions of per-hole scoring gains from these substitutions.

\begin{figure}[h!]
\centering
\includegraphics[width=0.6\textwidth]{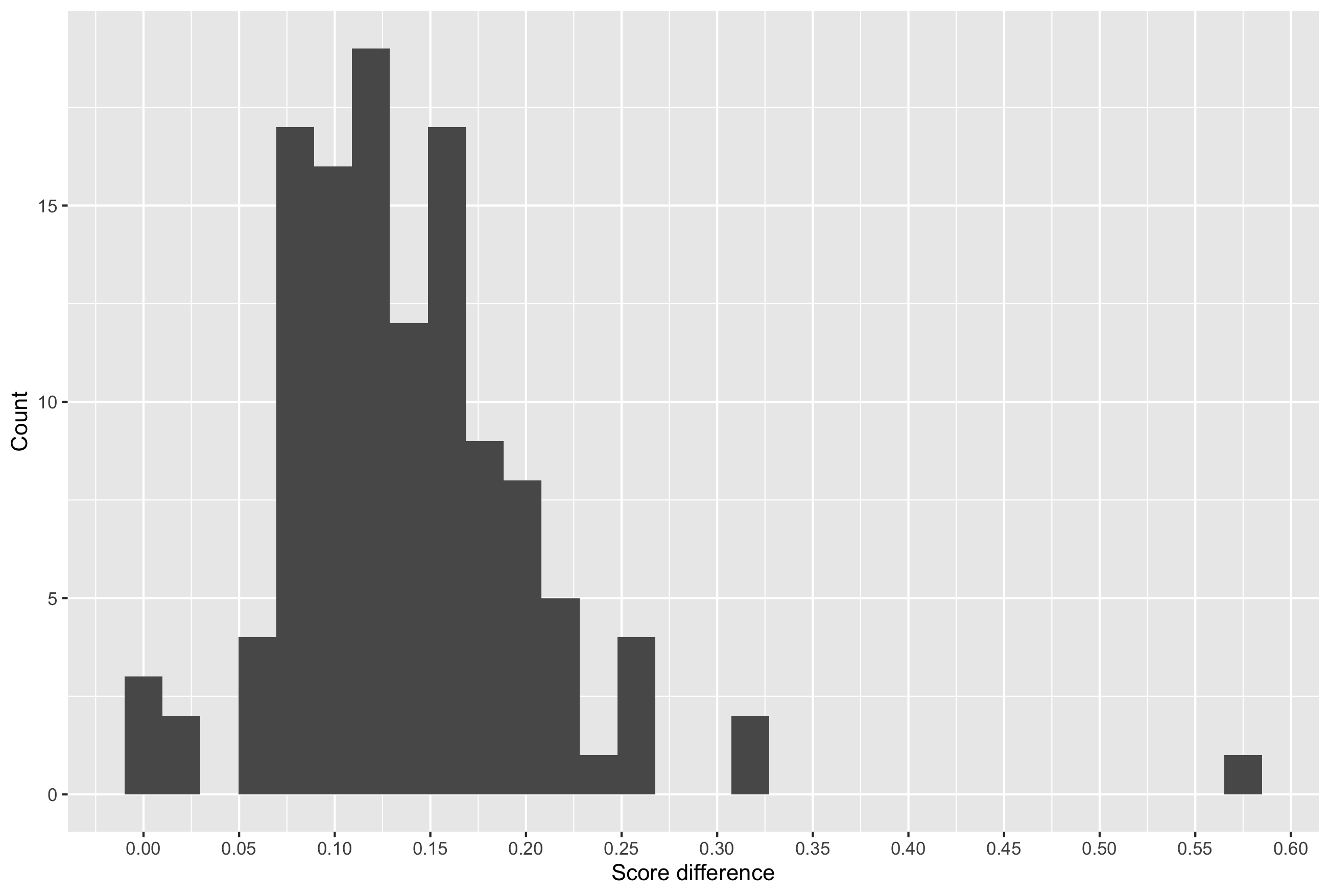}
\caption{Distribution of gain per hole (across 107 players) from substituting driving skills with McIlroy’s.}
\label{gap_driving}
\end{figure}

\begin{figure}[h!]
\centering
\includegraphics[width=0.6\textwidth]{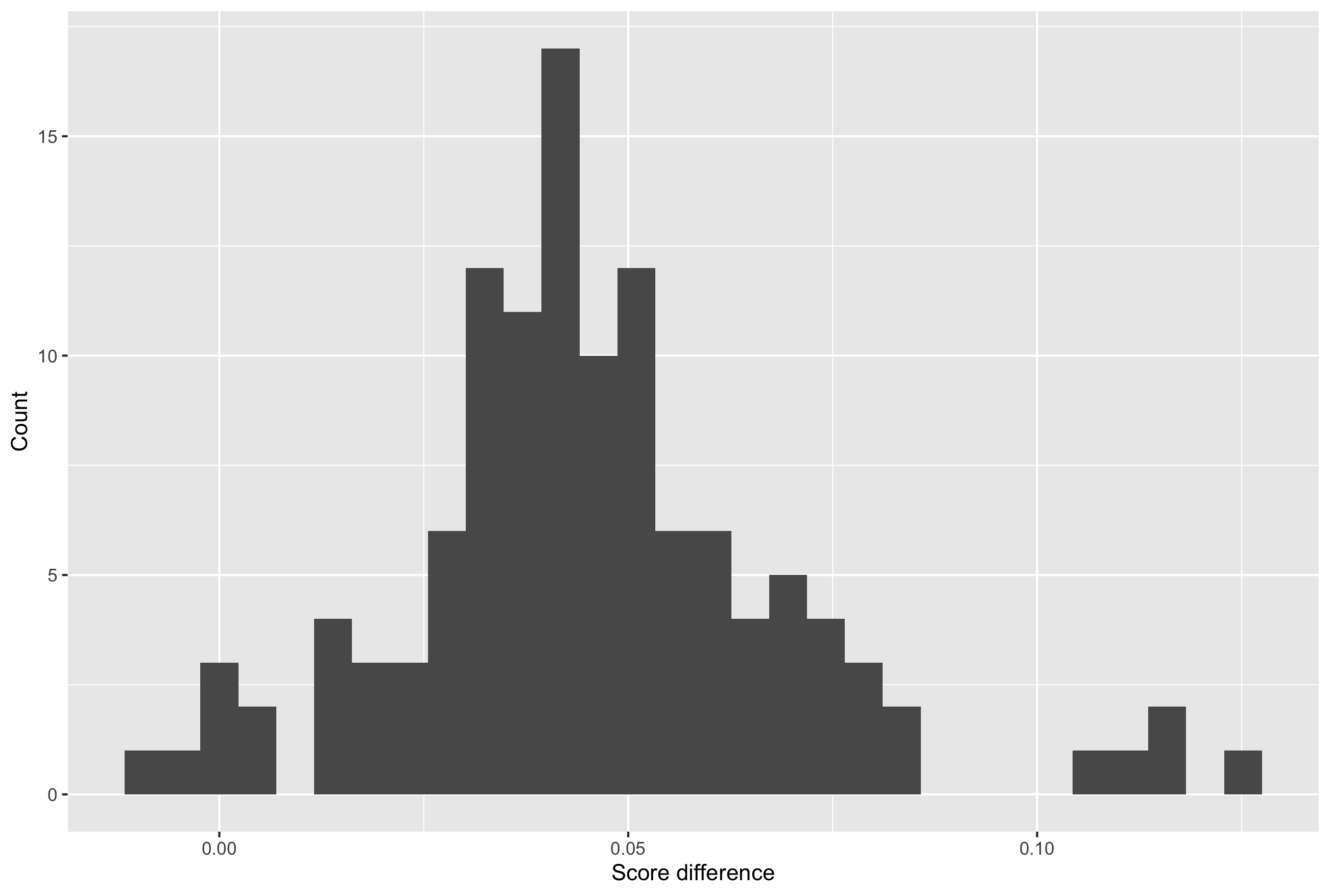}
\caption{Distribution of gain per hole (across 107 players) from substituting putting skills with Woods’s.}
\label{gap_putting}
\end{figure}

On average, adopting McIlroy’s driving yields a gain of 0.139 strokes per hole (95\% CI: [0.126, 0.152]), while adopting Woods’s putting results in a smaller gain of 0.046 strokes per hole (95\% CI: [0.041, 0.050]). These results suggest that, on average, gains in driving performance yield greater scoring benefits than improvements in putting (for the virtual players). This challenges the conventional wisdom encapsulated in the adage and aligns with previous empirical analyses questioning its validity~\cite{alexander2005drive,Baugher2016}.

Can these conclusions be extended to real PGA Tour players? We acknowledge that our model includes simplifications that prevent perfect replication of real-world performance. Tables~\ref{tab1} and~\ref{tab2} compare historical and simulated performance metrics for the top ten players at the 2018 Arnold Palmer Invitational (see Appendix Tables~\ref{tab3}--\ref{tab10} for 95\% confidence intervals). Metrics include total score, average tee-shot distance, fairway accuracy, directional misses, greens in regulation, and frequencies of water and bunker shots.

% latex table generated in R 4.2.2 by xtable 1.8-4 package
% Fri May  3 17:48:56 2024
\begin{table}[ht]
\centering
\begin{tabular}{rllrrrrrrrr}
  \hline
 & First Name & Last Name & Score & Tee-shot & Fairway & L & R & GiR  & Water & Bunker \\ 
  \hline
1 & Rory & McIlroy & 269 & 266.504 & 0.696 & 0.107 & 0.196 & 0.736 & 0.000 & 0.153 \\ 
  2 & Bryson & DeChambeau & 273 & 260.902 & 0.768 & 0.125 & 0.107 & 0.764 & 0.000 & 0.125 \\ 
  3 & Justin & Rose & 274 & 263.270 & 0.786 & 0.054 & 0.161 & 0.750 & 0.028 & 0.139 \\ 
  4 & Henrik & Stenson & 275 & 253.695 & 0.839 & 0.054 & 0.107 & 0.806 & 0.000 & 0.139 \\ 
  5 & Tiger & Woods & 277 & 258.218 & 0.643 & 0.107 & 0.250 & 0.708 & 0.000 & 0.181 \\ 
  6 & Ryan & Moore & 278 & 254.739 & 0.839 & 0.018 & 0.143 & 0.806 & 0.028 & 0.250 \\ 
  7 & Kevin & Chappell & 280 & 266.802 & 0.732 & 0.125 & 0.143 & 0.806 & 0.042 & 0.153 \\ 
  8 & Marc & Leishman & 280 & 259.035 & 0.714 & 0.071 & 0.214 & 0.722 & 0.014 & 0.208 \\ 
  9 & Patrick & Rodgers & 280 & 257.997 & 0.589 & 0.161 & 0.250 & 0.722 & 0.014 & 0.250 \\ 
  10 & Chris & Kirk & 281 & 254.156 & 0.732 & 0.143 & 0.125 & 0.708 & 0.014 & 0.264 \\ 
   \hline
\end{tabular}
  \caption{Historical golf metrics for the top ten players from the 2018 Arnold Palmer Invitational. Score: total score over the 4 rounds; Tee-shot : average tee-shot distance in meters (on par 4 and par 5 only); Fairway: percentage of fairways hit with the tee-shot on par 4 and par 5; L: percentage of fairways missed on the left on par 4 and par 5; R: percentage of fairways missed on the right on par 4 and par 5; GiR: percentage of green hit in a number of shot no more than par minus 2; Water: percentage of water hasard penalties ; Bunker:  percentage of bunker shots.}\label{tab1}
\end{table}

% latex table generated in R 4.2.2 by xtable 1.8-4 package
% Fri May  3 17:50:48 2024
\begin{table}[ht]
\centering
\begin{tabular}{rllrrrrrrrr}
  \hline
 & First Name & Last Name & Score & Tee-shot & Fairway & L & R & GiR  & Water & Bunker \\ 
  \hline
1 & Rory & McIlroy & 273.486 & 284.492 & 0.685 & 0.126 & 0.189 & 0.735 & 0.009 & 0.137 \\ 
  2 & Bryson & DeChambeau & 277.623 & 270.841 & 0.676 & 0.141 & 0.183 & 0.760 & 0.009 & 0.121 \\ 
  3 & Justin & Rose & 275.630 & 269.931 & 0.700 & 0.147 & 0.153 & 0.742 & 0.004 & 0.144 \\ 
  4 & Henrik & Stenson & 281.944 & 258.893 & 0.699 & 0.135 & 0.166 & 0.718 & 0.007 & 0.144 \\ 
  5 & Tiger & Woods & 269.085 & 263.598 & 0.765 & 0.088 & 0.147 & 0.795 & 0.005 & 0.099 \\ 
  6 & Ryan & Moore & 282.356 & 262.925 & 0.686 & 0.133 & 0.181 & 0.718 & 0.008 & 0.157 \\ 
  7 & Kevin & Chappell & 278.488 & 291.782 & 0.624 & 0.121 & 0.255 & 0.721 & 0.015 & 0.118 \\ 
  8 & Marc & Leishman & 274.260 & 277.578 & 0.702 & 0.136 & 0.162 & 0.770 & 0.009 & 0.101 \\ 
  9 & Patrick & Rodgers & 284.414 & 271.482 & 0.628 & 0.181 & 0.191 & 0.678 & 0.006 & 0.164 \\ 
  10 & Chris & Kirk & 289.085 & 243.917 & 0.647 & 0.135 & 0.217 & 0.693 & 0.013 & 0.119 \\ 
   \hline
\end{tabular}
\caption{Simulated golf metrics for the top ten players from the 2018 Arnold Palmer Invitational (10000 simulations). Score: average score over the 4 rounds. Tee-shot : average tee-shot distance in meters (on par 4 and par 5 only); Fairway: percentage of fairways hit with the tee-shot on par 4 and par 5; L: percentage of fairways missed on the left on par 4 and par 5; R: percentage of fairways missed on the right on par 4 and par 5; GiR: percentage of green hit in a number of shot no more than par minus 2; Water: percentage of water hasard penalties ; Bunker:  percentage of bunker shots }\label{tab2}
\end{table}

Several discrepancies emerge, most notably the model’s tendency to overestimate tee-shot distances. This is likely due to model simplifications—particularly its assumption of risk-neutral decision-making—rather than an overestimation of the players’ physical capabilities. In fact, the players are indeed capable of reaching the simulated distances. However, unlike real players, who may adopt more conservative strategies, the model consistently selects shots that minimize expected score. Targeting a longer landing zone increases the likelihood of missing the fairway and ending up in the rough, where outcomes are more variable. Although the model adjusts shot dispersion based on surface type and distance, it may still favor longer, riskier targets because dispersion generally decreases with shorter remaining distances, regardless of surface. As a result, the model may prefer a longer drive into the rough over a shorter, safer shot to the fairway. In contrast, a risk-averse player would adapt their strategy—shaping not only individual shot choices but also the overall approach from tee to pin—which helps explain discrepancies observed in other performance metrics.

Despite these limitations, the behaviors and skill profiles of our virtual players remain generally consistent with professional standards, as they are grounded in real data— TrackMan profiles derived from ShotLink—and informed by strategic assumptions reflective of PGA Tour play. To assess the model’s predictive accuracy, we evaluate how well it forecasts the average score per hole in the corresponding tournament for a cohort of PGA Tour players. Specifically, we randomly select 70 players to train a linear model that links simulated and observed scores, allowing for bias correction.

\begin{figure}[h!]
\centering
\includegraphics[width=0.6\textwidth]{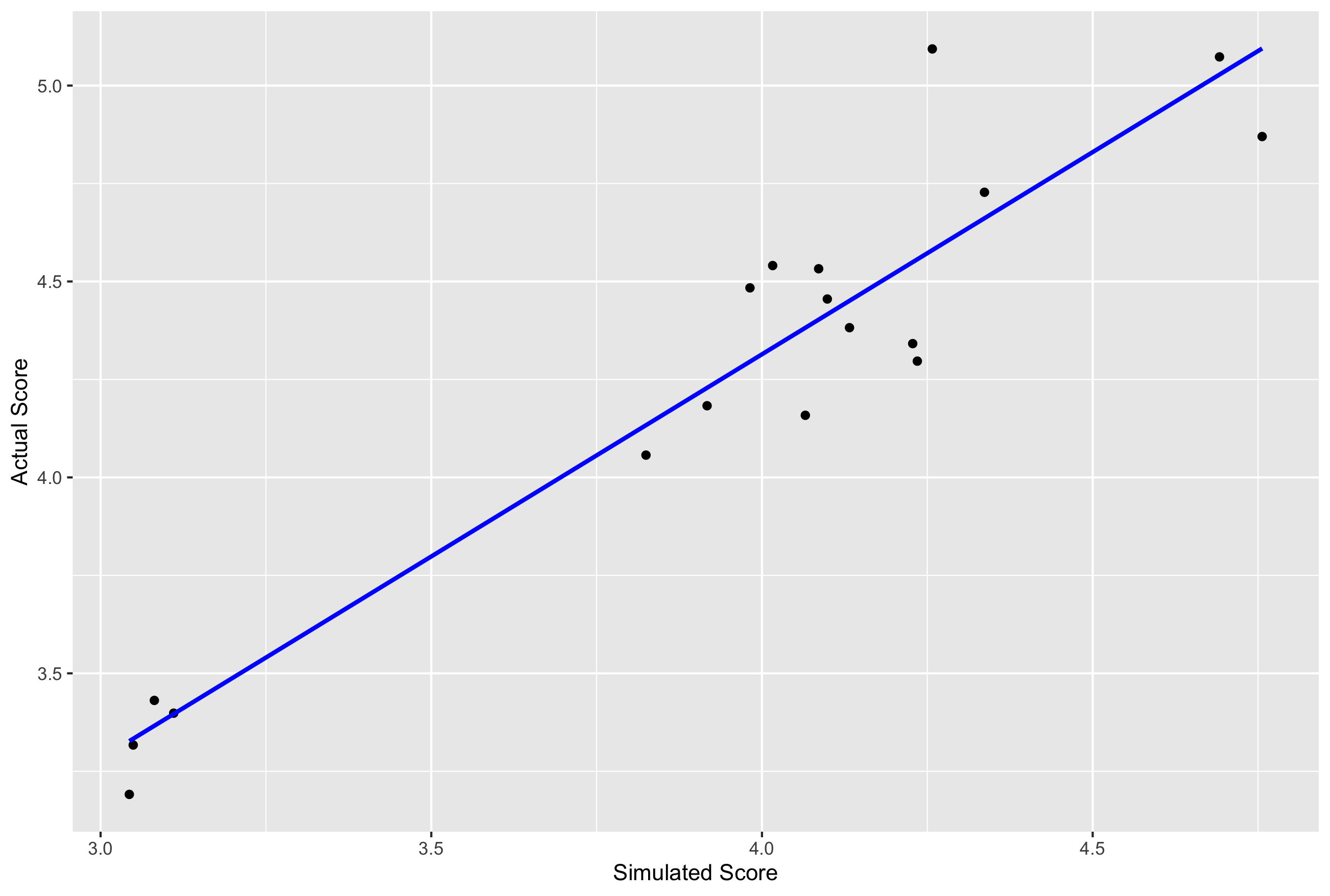}
\caption{Linear regression for predicting actual score for each of the 18 holes from simulated score (training set: 70 players).}
\label{linear_model}
\end{figure}

The in-sample $R^2$ of the model is 0.8927, with a slope of 1.03195 and an intercept of 0.18642, indicating a slight tendency to overestimate player performance. The out-of-sample $R^2$, computed on the remaining 37 players, is 0.8397. We use this model to forecast the average score of test cohorts and assess the variability of out-of-sample prediction error through repeated sampling. After 1,000 iterations, we obtain the distribution of prediction errors shown in Figure~\ref{forecast_error}, based on cohorts of 37 players.

\begin{figure}[h!]
\centering
\includegraphics[width=0.6\textwidth]{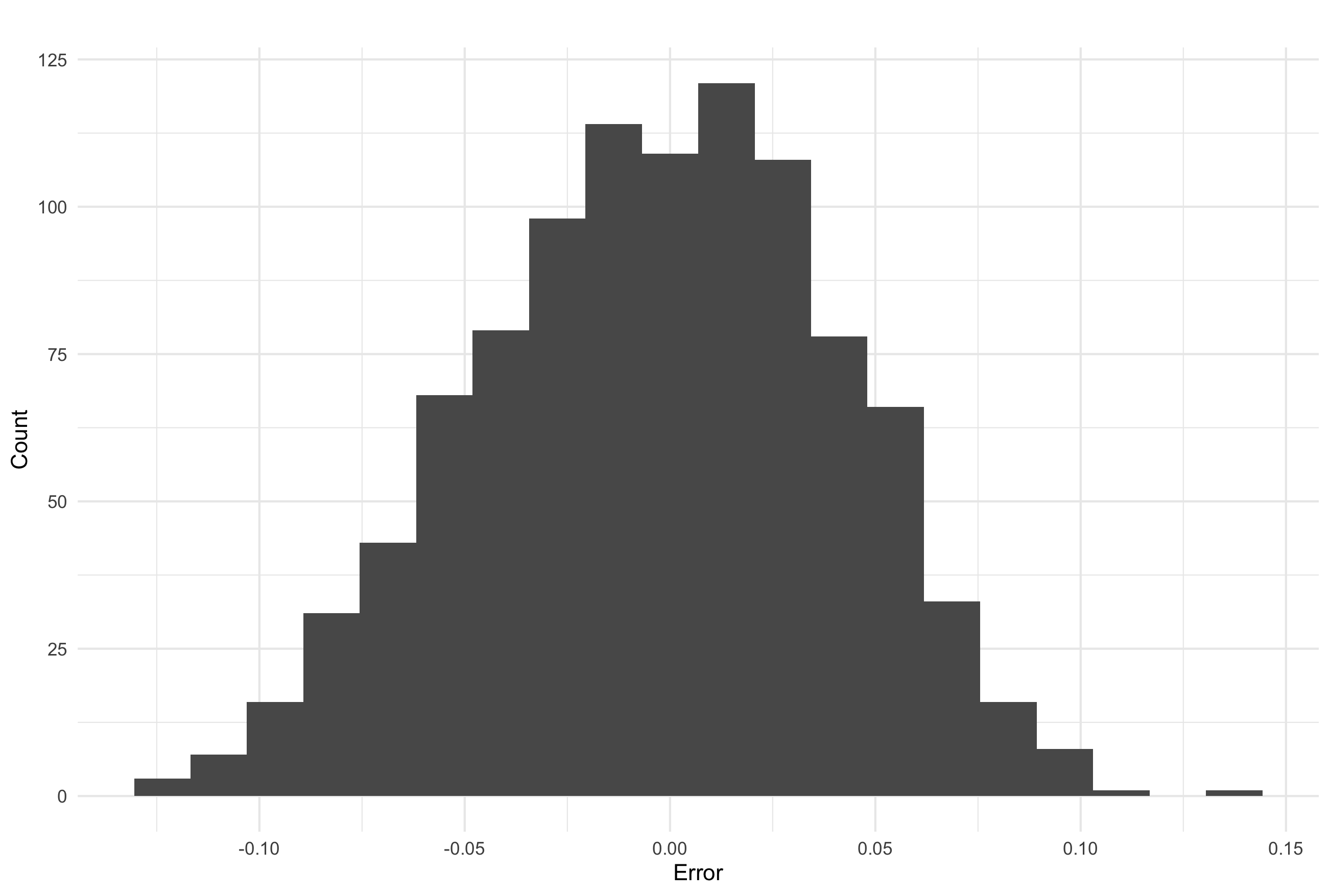}
\caption{Histogram of the out-of-sample error between forecasted and actual average score per hole for cohorts of 37 players.}
\label{forecast_error}
\end{figure}

The standard deviation of the error on the average score is $\sigma = 0.0435$. For a full cohort of 107 players, this scales to $\sqrt{37/107} \cdot \sigma \approx 0.0256$, yielding a 95\% confidence interval of approximately $\pm 0.05$. This confirms that the observed difference of $0.139 - 0.046 = 0.093$ strokes per hole—between adopting McIlroy’s driving and Woods’s putting for virtual players—which translates to an estimated difference of $(0.139 - 0.046) \times 1.03195 = 0.096$ strokes per hole for real players, is statistically significant. These results challenge the ``Drive for show, putt for dough'' adage not only in the virtual setting but also for actual PGA Tour players, with a high degree of confidence.

To conclude this section, we emphasize that, more broadly, this methodology—enhanced with a more detailed physics model—could be extended to the individual level and applied across different courses. It would enable the simulation of incremental improvements in key skills—such as reducing lateral dispersion, adding 10 meters to average driving distance, or increasing iron accuracy—and the assessment of their relative impact on overall performance. Such insights would support more informed, deliberate, and targeted investments in training and practice, not only for professionals but also for high-level amateurs seeking to optimize their development.

\section{Conclusion and Perspectives}\label{concl}

The primary aim of this study was to demonstrate the computational feasibility of an exact, data-driven methodology for optimizing strategic decision-making in professional golf. Our Markov Decision Process (MDP) framework, applied to PGA Tour ShotLink data, offers the first scalable implementation of a stochastic shortest path model in this context, highlighting that such large-scale models can be solved exactly with appropriate structuring and low-level optimization.

While predictive accuracy was not our primary objective, preliminary analyses suggest that the model yields relevant and interpretable outcomes. For example, we have used the model to explore the marginal value of specific skills by performing controlled substitutions—replacing players' driving or putting abilities with those of elite performers. This allowed us to critically examine the popular adage ``Drive for show, putt for dough.'' Together, these experiments offer initial validation that the model captures meaningful behavioral patterns and, with further refinement, could be a valuable tool for skill assessment and training prioritization. In addition, our virtual players could also be used to evaluate course design proposals, simulate the impact of structural changes (e.g., repositioning hazards or tees), or rank courses based on simulated performance metrics for a representative cohort. 

We emphasize, however, that this work remains a proof of concept. The structure of the model is general and flexible, but its current implementation involves several simplifying assumptions: straight-line trajectories, no ball roll, infinite tree height, homogeneous and flat green surfaces, and “average” weather conditions. While such assumptions are consistent with the current state of the art, it is clear that precision could benefit from a more sophisticated physics engine. Yet implementing such refinements would involve substantial additional engineering effort—already, efficiently simulating the current model posed significant challenges.

Importantly, these physical enhancements pose no conceptual difficulty within our MDP framework. They can be accommodated by adjusting the model’s state, action, and transition components, though maintaining tractability would require additional computational strategies.  Potential avenues include extending the Bresenham algorithm to 3D \cite{AU}, exploiting sparsity in the transition structure (e.g., since each action can only reach a limited neighborhood), and applying low-level optimization techniques such as parallelization and memory caching. Although such techniques were not implemented in the current study, we are confident that—with careful engineering—the computational demands of more realistic models could be kept under control.

Although our empirical study is based on professional golf data, the methodology is fully compatible with amateur golfers, provided that sufficiently granular tracking data is available—such as that captured by systems like Arccos or TrackMan. The model could then be used to simulate expected performance under optimal strategies personalized to the player's own historical data, offering meaningful insights for recreational players and coaches alike. Furthermore, the framework is adaptable to other competition formats, such as match play, as explored in a separate study \cite{wajge2024should}.

Finally, one may argue that extensions of the model should incorporate behavioral dynamics such as momentum or the so-called ``hot hand'' phenomenon. In this study, we adopt the Markov property—treating each action as independent—which is consistent with both standard practice in sequential decision modeling and the professional mindset of evaluating each shot in isolation (see, for example, \cite{rotella2012golf}). Nonetheless, we acknowledge that positive or negative performance trends may influence outcomes in practice. While empirical evidence for the ``hot hand'' effect in golf is mixed, studies suggest that the reverse effect—the ``cold hand''—may be more robustly supported by data \cite{elmore2018hot}. From a modeling perspective, incorporating such dynamics would be conceptually straightforward, for instance by extending the state space to capture recent shot history. However, doing so would introduce additional complexity and would require richer empirical data to calibrate and validate the behavior of virtual players accordingly.

In summary, we provide a scalable, transparent, and reproducible framework for golf strategy modeling, with open-source code in R and C++. This work lays the foundation for future research that may integrate finer physical models, richer behavioral dynamics, or advanced validation methods. We hope it will contribute to both academic exploration and practical decision-making in the evolving field of sports analytics.

\section{Acknowledgement}

We would also like to thank Renaud Gris and Jason Belot (French U21 Elite amateur trainers) for the constructive discussions of our models and providing data on some of their elite amateurs. We would like to thank the PGA Tour for giving us access to the ShotLink data. At the time of writing, the ShotLink Intelligence program has been discontinued, and ShotLink data is no longer publicly available to academics. However, we have received authorization from the PGA Tour to make our data accessible for replication and validation by other researchers. We would like to express our gratitude to the PGA Tour, and to Ken Lovell in particular, for their support. Please note that any further use of the data requires prior permission from the PGA Tour. We are very grateful to Marc-Olivier Boldi for the insightful exchanges on how to statistically demonstrate that our conclusions for virtual players could be extended to real players through cohort-based analysis.

\bibliographystyle{plain} % You can choose a different bibliography style if needed
\bibliography{paper}  % Use the name of your .bib file here (without the .bib extension)

\newpage \section*{Appendix}

% latex table generated in R 4.2.2 by xtable 1.8-4 package
% Wed Sep 11 16:52:38 2024
\begin{table}[!htbp]
\centering
\begin{tabular}{rllrrrrrr}
  \hline
 & FirstName & LastName & Score & Score LB & score UB & Tee-shot & Tee-shot LB & Tee-shot UB \\ 
  \hline
1 & Rory & McIlroy & 269.000 & 258.816 & 279.184 & 266.504 & 262.347 & 270.661 \\ 
  2 & Bryson & DeChambeau & 273.000 & 262.445 & 283.555 & 260.902 & 257.769 & 264.036 \\ 
  3 & Justin & Rose & 274.000 & 263.146 & 284.854 & 263.270 & 259.396 & 267.145 \\ 
  4 & Henrik & Stenson & 275.000 & 266.162 & 283.838 & 253.695 & 249.943 & 257.446 \\ 
  5 & Tiger & Woods & 277.000 & 267.877 & 286.123 & 258.218 & 253.077 & 263.359 \\ 
  6 & Ryan & Moore & 278.000 & 268.006 & 287.994 & 254.739 & 250.969 & 258.510 \\ 
  7 & Kevin & Chappell & 280.000 & 268.684 & 291.316 & 266.802 & 262.215 & 271.389 \\ 
  8 & Marc & Leishman & 280.000 & 269.146 & 290.854 & 259.035 & 254.889 & 263.180 \\ 
  9 & Patrick & Rodgers & 280.000 & 270.532 & 289.468 & 257.997 & 254.502 & 261.492 \\ 
  10 & Chris & Kirk & 281.000 & 270.205 & 291.795 & 254.156 & 250.510 & 257.801 \\ 
%  11 & Rickie & Fowler & 283.000 & 272.942 & 293.058 & 259.907 & 256.295 & 263.518 \\ 
%  12 & Hideki & Matsuyama & 287.000 & 276.445 & 297.555 & 263.563 & 258.687 & 268.439 \\ 
%  13 & Sam & Burns & 288.000 & 274.901 & 301.099 & 266.190 & 263.166 & 269.213 \\ 
%  14 & Kiradech & Aphibarnrat & 292.000 & 280.913 & 303.087 & 251.284 & 248.174 & 254.394 \\ 
%  15 & Matthew & Fitzpatrick & 294.000 & 280.141 & 307.859 & 243.170 & 236.303 & 250.038 \\ 
%  16 & Sangmoon & Bae & 298.000 & 287.265 & 308.735 & 254.137 & 250.525 & 257.750 \\ 
   \hline
\end{tabular}
\caption{The table presents  the confidence intervals for the mean of the Score and Tee-shot metrics considered in Table  \ref{tab1}. For each metric, 'LB' represents the lower bound and 'UB' represents the upper bound of the interval.}\label{tab3}
\end{table}

% latex table generated in R 4.2.2 by xtable 1.8-4 package
% Wed Sep 11 16:52:38 2024
\begin{table}[!htbp]
\centering
\begin{tabular}{rllrrrrrr}
  \hline
 & FirstName & LastName & Fairway & Fairway LB & Fairway UB & L & L LB & L UB \\ 
  \hline
1 & Rory & McIlroy & 0.696 & 0.594 & 0.799 & 0.107 & 0.036 & 0.178 \\ 
  2 & Bryson & DeChambeau & 0.768 & 0.675 & 0.860 & 0.125 & 0.061 & 0.189 \\ 
  3 & Justin & Rose & 0.786 & 0.691 & 0.880 & 0.054 & 0.000 & 0.107 \\ 
  4 & Henrik & Stenson & 0.839 & 0.758 & 0.921 & 0.054 & 0.006 & 0.101 \\ 
  5 & Tiger & Woods & 0.643 & 0.536 & 0.750 & 0.107 & 0.036 & 0.178 \\ 
  6 & Ryan & Moore & 0.839 & 0.754 & 0.925 & 0.018 & -0.013 & 0.049 \\ 
  7 & Kevin & Chappell & 0.732 & 0.636 & 0.828 & 0.125 & 0.047 & 0.203 \\ 
  8 & Marc & Leishman & 0.714 & 0.620 & 0.809 & 0.071 & 0.015 & 0.128 \\ 
  9 & Patrick & Rodgers & 0.589 & 0.500 & 0.678 & 0.161 & 0.087 & 0.234 \\ 
  10 & Chris & Kirk & 0.732 & 0.640 & 0.825 & 0.143 & 0.067 & 0.218 \\ 
%  11 & Rickie & Fowler & 0.643 & 0.539 & 0.747 & 0.143 & 0.067 & 0.218 \\ 
%  12 & Hideki & Matsuyama & 0.607 & 0.489 & 0.725 & 0.250 & 0.149 & 0.351 \\ 
%  13 & Sam & Burns & 0.714 & 0.613 & 0.815 & 0.107 & 0.032 & 0.183 \\ 
%  14 & Kiradech & Aphibarnrat & 0.679 & 0.563 & 0.794 & 0.071 & 0.010 & 0.133 \\ 
%  15 & Matthew & Fitzpatrick & 0.643 & 0.536 & 0.750 & 0.179 & 0.081 & 0.276 \\ 
%  16 & Sangmoon & Bae & 0.536 & 0.460 & 0.611 & 0.250 & 0.174 & 0.326 \\ 
   \hline
\end{tabular}
\caption{The table presents  the confidence intervals for the mean of the Fairway and L metrics considered in Table  \ref{tab1}. For each metric, 'LB' represents the lower bound and 'UB' represents the upper bound of the interval.}\label{tab4}

\end{table}

% latex table generated in R 4.2.2 by xtable 1.8-4 package
% Wed Sep 11 16:52:38 2024
\begin{table}[!htbp]
\centering
\begin{tabular}{rllrrrrrr}
  \hline
& FirstName & LastName & R & R LB & R UB & GiR & GiR LB & GiR UB \\ 
  \hline
1 & Rory & McIlroy & 0.196 & 0.107 & 0.286 & 0.736 & 0.638 & 0.834 \\ 
  2 & Bryson & DeChambeau & 0.107 & 0.036 & 0.178 & 0.764 & 0.682 & 0.846 \\ 
  3 & Justin & Rose & 0.161 & 0.075 & 0.246 & 0.750 & 0.648 & 0.852 \\ 
  4 & Henrik & Stenson & 0.107 & 0.032 & 0.183 & 0.806 & 0.725 & 0.886 \\ 
  5 & Tiger & Woods & 0.250 & 0.146 & 0.354 & 0.708 & 0.610 & 0.806 \\ 
  6 & Ryan & Moore & 0.143 & 0.063 & 0.223 & 0.806 & 0.729 & 0.883 \\ 
  7 & Kevin & Chappell & 0.143 & 0.072 & 0.214 & 0.806 & 0.714 & 0.897 \\ 
  8 & Marc & Leishman & 0.214 & 0.123 & 0.305 & 0.722 & 0.618 & 0.826 \\ 
  9 & Patrick & Rodgers & 0.250 & 0.159 & 0.341 & 0.722 & 0.623 & 0.822 \\ 
  10 & Chris & Kirk & 0.125 & 0.052 & 0.198 & 0.708 & 0.605 & 0.811 \\ 
%  11 & Rickie & Fowler & 0.214 & 0.123 & 0.305 & 0.778 & 0.683 & 0.872 \\ 
%  12 & Hideki & Matsuyama & 0.143 & 0.059 & 0.226 & 0.597 & 0.476 & 0.718 \\ 
%  13 & Sam & Burns & 0.179 & 0.107 & 0.250 & 0.681 & 0.564 & 0.797 \\ 
%  14 & Kiradech & Aphibarnrat & 0.250 & 0.152 & 0.348 & 0.583 & 0.456 & 0.711 \\ 
%  15 & Matthew & Fitzpatrick & 0.179 & 0.081 & 0.276 & 0.639 & 0.537 & 0.741 \\ 
%  16 & Sangmoon & Bae & 0.214 & 0.127 & 0.302 & 0.611 & 0.517 & 0.705 \\ 
   \hline
\end{tabular}
\caption{The table presents  the confidence intervals for the mean of the R and GiR metrics considered in Table  \ref{tab1}. For each metric, 'LB' represents the lower bound and 'UB' represents the upper bound of the interval.}\label{tab5}

\end{table}

% latex table generated in R 4.2.2 by xtable 1.8-4 package
% Wed Sep 11 16:52:38 2024
\begin{table}[!htbp]
\centering
\begin{tabular}{rllrrrrrr}
  \hline
 & FirstName & LastName & water & water LB & water UB & bunker & bunker LB & bunker UB \\ 
  \hline
1 & Rory & McIlroy & 0.000 & 0.000 & 0.000 & 0.153 & 0.065 & 0.240 \\ 
  2 & Bryson & DeChambeau & 0.000 & 0.000 & 0.000 & 0.125 & 0.046 & 0.204 \\ 
  3 & Justin & Rose & 0.028 & -0.004 & 0.059 & 0.139 & 0.056 & 0.222 \\ 
  4 & Henrik & Stenson & 0.000 & 0.000 & 0.000 & 0.139 & 0.069 & 0.209 \\ 
  5 & Tiger & Woods & 0.000 & 0.000 & 0.000 & 0.181 & 0.093 & 0.268 \\ 
  6 & Ryan & Moore & 0.028 & -0.027 & 0.082 & 0.250 & 0.151 & 0.349 \\ 
  7 & Kevin & Chappell & 0.042 & -0.019 & 0.103 & 0.153 & 0.062 & 0.243 \\ 
  8 & Marc & Leishman & 0.014 & -0.013 & 0.041 & 0.208 & 0.115 & 0.301 \\ 
  9 & Patrick & Rodgers & 0.014 & -0.013 & 0.041 & 0.250 & 0.158 & 0.342 \\ 
  10 & Chris & Kirk & 0.014 & -0.013 & 0.041 & 0.264 & 0.161 & 0.367 \\ 
%  11 & Rickie & Fowler & 0.042 & -0.005 & 0.089 & 0.139 & 0.065 & 0.213 \\ 
%  12 & Hideki & Matsuyama & 0.028 & -0.011 & 0.066 & 0.236 & 0.131 & 0.342 \\ 
%  13 & Sam & Burns & 0.042 & -0.005 & 0.089 & 0.292 & 0.163 & 0.420 \\ 
%  14 & Kiradech & Aphibarnrat & 0.028 & -0.011 & 0.066 & 0.278 & 0.156 & 0.400 \\ 
%  15 & Matthew & Fitzpatrick & 0.083 & 0.017 & 0.150 & 0.278 & 0.144 & 0.411 \\ 
%  16 & Sangmoon & Bae & 0.056 & 0.001 & 0.110 & 0.389 & 0.245 & 0.533 \\ 
   \hline
\end{tabular}
\caption{The table presents  the confidence intervals for the mean of the water and bunker metrics considered in Table  \ref{tab1}. For each metric, 'LB' represents the lower bound and 'UB' represents the upper bound of the interval.}\label{tab6}

\end{table}

% latex table generated in R 4.2.2 by xtable 1.8-4 package
% Wed Sep 11 16:49:02 2024
\begin{table}[!htbp]
\centering
\begin{tabular}{rllrrrrrr}
  \hline
 & FirstName & LastName & Score & Score LB & score UB & Tee-shot & Tee-shot LB & Tee-shot UB \\ 
  \hline
1 & Rory & McIlroy & 273.486 & 272.555 & 274.418 & 284.492 & 283.870 & 285.113 \\ 
  2 & Bryson & DeChambeau & 277.623 & 276.739 & 278.506 & 270.841 & 270.328 & 271.354 \\ 
  3 & Justin & Rose & 275.630 & 274.768 & 276.493 & 269.931 & 269.510 & 270.353 \\ 
  4 & Henrik & Stenson & 281.944 & 281.022 & 282.867 & 258.893 & 258.463 & 259.323 \\ 
  5 & Tiger & Woods & 269.085 & 268.275 & 269.895 & 263.598 & 263.269 & 263.927 \\ 
  6 & Ryan & Moore & 282.356 & 281.473 & 283.240 & 262.925 & 262.520 & 263.330 \\ 
  7 & Kevin & Chappell & 278.488 & 277.572 & 279.405 & 291.782 & 291.217 & 292.348 \\ 
  8 & Marc & Leishman & 274.260 & 273.386 & 275.135 & 277.578 & 277.045 & 278.110 \\ 
  9 & Patrick & Rodgers & 284.414 & 283.510 & 285.317 & 271.482 & 271.055 & 271.908 \\ 
  10 & Chris & Kirk & 289.085 & 288.156 & 290.013 & 243.917 & 243.217 & 244.617 \\ 
%  11 & Rickie & Fowler & 274.224 & 273.374 & 275.075 & 265.435 & 265.085 & 265.785 \\ 
%  12 & Hideki & Matsuyama & 275.598 & 274.754 & 276.442 & 270.190 & 269.780 & 270.600 \\ 
%  13 & Sam & Burns & 274.419 & 273.555 & 275.283 & 258.712 & 258.320 & 259.104 \\ 
%  14 & Kiradech & Aphibarnrat & 273.199 & 272.343 & 274.054 & 272.418 & 272.198 & 272.637 \\ 
%  15 & Matthew & Fitzpatrick & 274.724 & 273.836 & 275.613 & 257.578 & 257.246 & 257.911 \\ 
%  16 & Sangmoon & Bae & 270.990 & 270.169 & 271.811 & 271.870 & 271.471 & 272.269 \\ 
   \hline
\end{tabular}
\caption{The table presents  the confidence intervals for the mean of the Score and Tee-shot metrics considered in Table  \ref{tab2}. For each metric, 'LB' represents the lower bound and 'UB' represents the upper bound of the interval.}\label{tab7}
\end{table}

% latex table generated in R 4.2.2 by xtable 1.8-4 package
% Wed Sep 11 16:49:02 2024
\begin{table}[!htbp]
\centering
\begin{tabular}{rllrrrrrr}
  \hline
 & FirstName & LastName & Fairway & Fairway LB & Fairway UB & L & L LB & L UB \\ 
  \hline
1 & Rory & McIlroy & 0.685 & 0.677 & 0.693 & 0.126 & 0.120 & 0.133 \\ 
  2 & Bryson & DeChambeau & 0.676 & 0.668 & 0.685 & 0.141 & 0.135 & 0.147 \\ 
  3 & Justin & Rose & 0.700 & 0.691 & 0.709 & 0.147 & 0.140 & 0.154 \\ 
  4 & Henrik & Stenson & 0.699 & 0.690 & 0.707 & 0.135 & 0.129 & 0.141 \\ 
  5 & Tiger & Woods & 0.765 & 0.757 & 0.773 & 0.088 & 0.083 & 0.093 \\ 
  6 & Ryan & Moore & 0.686 & 0.677 & 0.695 & 0.133 & 0.126 & 0.139 \\ 
  7 & Kevin & Chappell & 0.624 & 0.615 & 0.632 & 0.121 & 0.115 & 0.127 \\ 
  8 & Marc & Leishman & 0.702 & 0.694 & 0.711 & 0.136 & 0.130 & 0.142 \\ 
  9 & Patrick & Rodgers & 0.628 & 0.619 & 0.637 & 0.181 & 0.174 & 0.188 \\ 
  10 & Chris & Kirk & 0.647 & 0.639 & 0.656 & 0.135 & 0.129 & 0.142 \\ 
%  11 & Rickie & Fowler & 0.686 & 0.677 & 0.694 & 0.132 & 0.126 & 0.138 \\ 
%  12 & Hideki & Matsuyama & 0.745 & 0.737 & 0.753 & 0.105 & 0.099 & 0.110 \\ 
%  13 & Sam & Burns & 0.705 & 0.697 & 0.713 & 0.146 & 0.140 & 0.153 \\ 
%  14 & Kiradech & Aphibarnrat & 0.829 & 0.822 & 0.836 & 0.068 & 0.064 & 0.073 \\ 
%  15 & Matthew & Fitzpatrick & 0.807 & 0.799 & 0.814 & 0.073 & 0.069 & 0.078 \\ 
%  16 & Sangmoon & Bae & 0.890 & 0.884 & 0.895 & 0.049 & 0.045 & 0.053 \\ 
   \hline
\end{tabular}
\caption{The table presents  the confidence intervals for the mean of the Fairway and L metrics considered in Table  \ref{tab2}. For each metric, 'LB' represents the lower bound and 'UB' represents the upper bound of the interval.}\label{tab8}
\end{table}
% latex table generated in R 4.2.2 by xtable 1.8-4 package
% Wed Sep 11 16:49:02 2024
\begin{table}[!htbp]
\centering
\begin{tabular}{rllrrrrrr}
  \hline
& FirstName & LastName & R & R LB & R UB & GiR & GiR LB & GiR UB \\ 
  \hline
1 & Rory & McIlroy & 0.189 & 0.182 & 0.196 & 0.735 & 0.726 & 0.743 \\ 
  2 & Bryson & DeChambeau & 0.183 & 0.176 & 0.190 & 0.760 & 0.752 & 0.768 \\ 
  3 & Justin & Rose & 0.153 & 0.146 & 0.160 & 0.742 & 0.734 & 0.750 \\ 
  4 & Henrik & Stenson & 0.166 & 0.159 & 0.173 & 0.718 & 0.710 & 0.727 \\ 
  5 & Tiger & Woods & 0.147 & 0.141 & 0.154 & 0.795 & 0.787 & 0.803 \\ 
  6 & Ryan & Moore & 0.181 & 0.175 & 0.188 & 0.718 & 0.710 & 0.727 \\ 
  7 & Kevin & Chappell & 0.255 & 0.248 & 0.262 & 0.721 & 0.713 & 0.730 \\ 
  8 & Marc & Leishman & 0.162 & 0.155 & 0.168 & 0.770 & 0.761 & 0.778 \\ 
  9 & Patrick & Rodgers & 0.191 & 0.184 & 0.198 & 0.678 & 0.669 & 0.687 \\ 
  10 & Chris & Kirk & 0.217 & 0.210 & 0.224 & 0.693 & 0.684 & 0.701 \\ 
%  11 & Rickie & Fowler & 0.183 & 0.176 & 0.190 & 0.773 & 0.765 & 0.781 \\ 
%  12 & Hideki & Matsuyama & 0.150 & 0.144 & 0.157 & 0.781 & 0.773 & 0.789 \\ 
%  13 & Sam & Burns & 0.149 & 0.142 & 0.155 & 0.764 & 0.756 & 0.772 \\ 
%  14 & Kiradech & Aphibarnrat & 0.103 & 0.098 & 0.109 & 0.762 & 0.753 & 0.770 \\ 
%  15 & Matthew & Fitzpatrick & 0.120 & 0.114 & 0.126 & 0.749 & 0.741 & 0.757 \\ 
%  16 & Sangmoon & Bae & 0.061 & 0.057 & 0.066 & 0.830 & 0.823 & 0.837 \\ 
   \hline
\end{tabular}
\caption{The table presents  the confidence intervals for the mean of the R and GiR metrics considered in Table  \ref{tab2}. For each metric, 'LB' represents the lower bound and 'UB' represents the upper bound of the interval.}\label{tab9}
\end{table}
% latex table generated in R 4.2.2 by xtable 1.8-4 package
% Wed Sep 11 16:49:02 2024
% latex table generated in R 4.2.2 by xtable 1.8-4 package
% Wed Sep 11 17:07:01 2024
\begin{table}[!htbp]
\centering
\begin{tabular}{rllrrrrrr}
  \hline
 & FirstName & LastName & water & water LB & water UB & bunker & bunker LB & bunker UB \\ 
  \hline
1 & Rory & McIlroy & 0.009 & 0.007 & 0.011 & 0.137 & 0.130 & 0.144 \\ 
  2 & Bryson & DeChambeau & 0.009 & 0.007 & 0.010 & 0.121 & 0.114 & 0.128 \\ 
  3 & Justin & Rose & 0.004 & 0.003 & 0.005 & 0.144 & 0.137 & 0.150 \\ 
  4 & Henrik & Stenson & 0.007 & 0.006 & 0.009 & 0.144 & 0.137 & 0.151 \\ 
  5 & Tiger & Woods & 0.005 & 0.003 & 0.006 & 0.099 & 0.093 & 0.105 \\ 
  6 & Ryan & Moore & 0.008 & 0.006 & 0.009 & 0.157 & 0.150 & 0.164 \\ 
  7 & Kevin & Chappell & 0.015 & 0.013 & 0.017 & 0.118 & 0.112 & 0.125 \\ 
  8 & Marc & Leishman & 0.009 & 0.007 & 0.011 & 0.101 & 0.095 & 0.107 \\ 
  9 & Patrick & Rodgers & 0.006 & 0.005 & 0.008 & 0.164 & 0.157 & 0.172 \\ 
  10 & Chris & Kirk & 0.013 & 0.011 & 0.015 & 0.107 & 0.113 & 0.126 \\ 
%  11 & Rickie & Fowler & 0.003 & 0.002 & 0.004 & 0.125 & 0.118 & 0.131 \\ 
%  12 & Hideki & Matsuyama & 0.002 & 0.001 & 0.003 & 0.110 & 0.103 & 0.116 \\ 
%  13 & Sam & Burns & 0.008 & 0.006 & 0.010 & 0.136 & 0.129 & 0.143 \\ 
%  14 & Kiradech & Aphibarnrat & 0.003 & 0.002 & 0.004 & 0.139 & 0.132 & 0.146 \\ 
%  15 & Matthew & Fitzpatrick & 0.004 & 0.003 & 0.006 & 0.132 & 0.125 & 0.139 \\ 
%  16 & Sangmoon & Bae & 0.003 & 0.002 & 0.004 & 0.087 & 0.081 & 0.092 \\ 
   \hline
\end{tabular}
\caption{The table presents  the confidence intervals for the mean of the water and bunker metrics considered in Table  \ref{tab2}. For each metric, 'LB' represents the lower bound and 'UB' represents the upper bound of the interval.}\label{tab10}
\end{table}

\begin{figure*}[!htbp]
\includegraphics[width=0.3\linewidth]{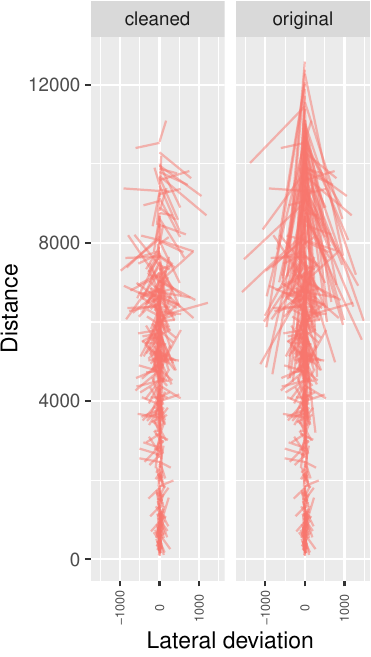} \caption{\label{fig:philrough} Each segment of the figures represents a target/destination pair for shots played from the rough. All target points have been rotated so as to appear on the y-axis. The ``original'' panel shows the inferred data for Phil Mickelson before removing outliers. The ``cleaned'' panel shows the data for Phil Mickelson after removal of outliers. Numbers are in meters.}\label{fig:unnamed-chunk-3}
\end{figure*}

\begin{figure*}[!htbp]
\includegraphics{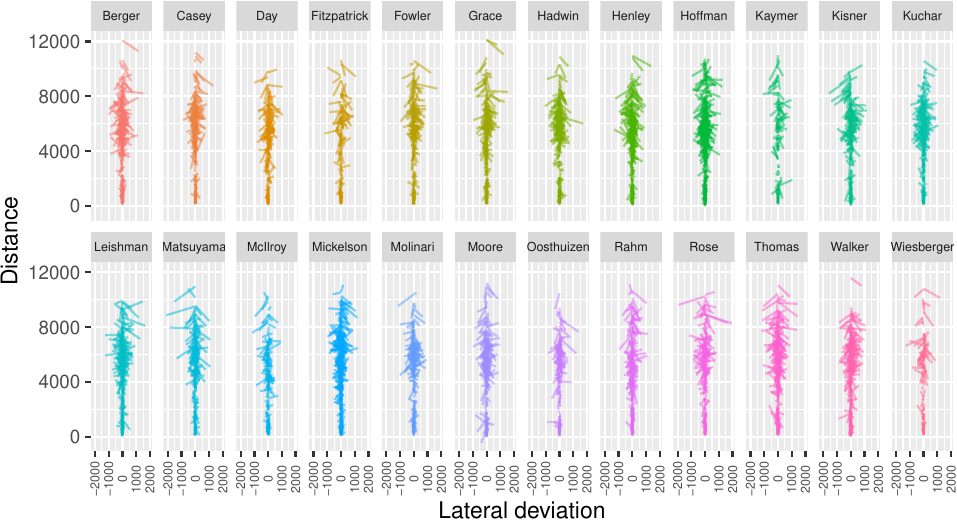} \caption{\label{fig:someplayer1rough} Each segment of the figure represents a target/destination pair for shots played from the rough after removal of outliers. All target points have been rotated so as to appear on the y-axis.  Numbers are in meters. }\label{fig:unnamed-chunk-4}
\end{figure*}

\begin{figure*}[!htbp]
\includegraphics[height=0.4\textheight]{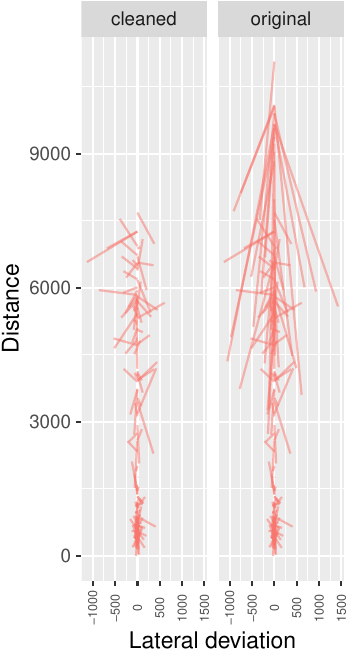} \caption{\label{fig:philbunker} Each segment on the figures represents a target/destination pair for shots played from bunkers. All target points have been rotated so as to appear on the y-axis. The ``original'' panel shows the inferred data for Phil Mickelson before removing outliers. The ``cleaned'' panel shows the data for Phil Mickelson after removal of outliers. Numbers are in meters.}\label{fig:unnamed-chunk-5}
\end{figure*}

\begin{figure*}[!htbp]
\includegraphics{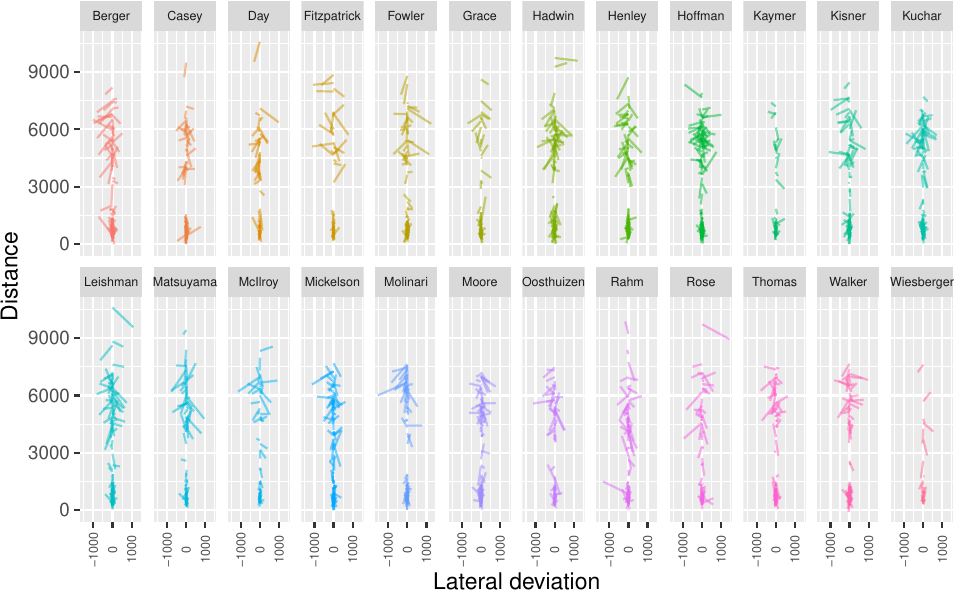} \caption{\label{fig:someplayer1bunker} Each segment of the figure represents a target/destination pair for shots played from the bunker after removal of outliers. All target points have been rotated so as to appear on the y-axis.  Numbers are in meters. }\label{fig:unnamed-chunk-6}
\end{figure*}

\begin{figure*}[!htbp]
\includegraphics[width=0.5\linewidth]{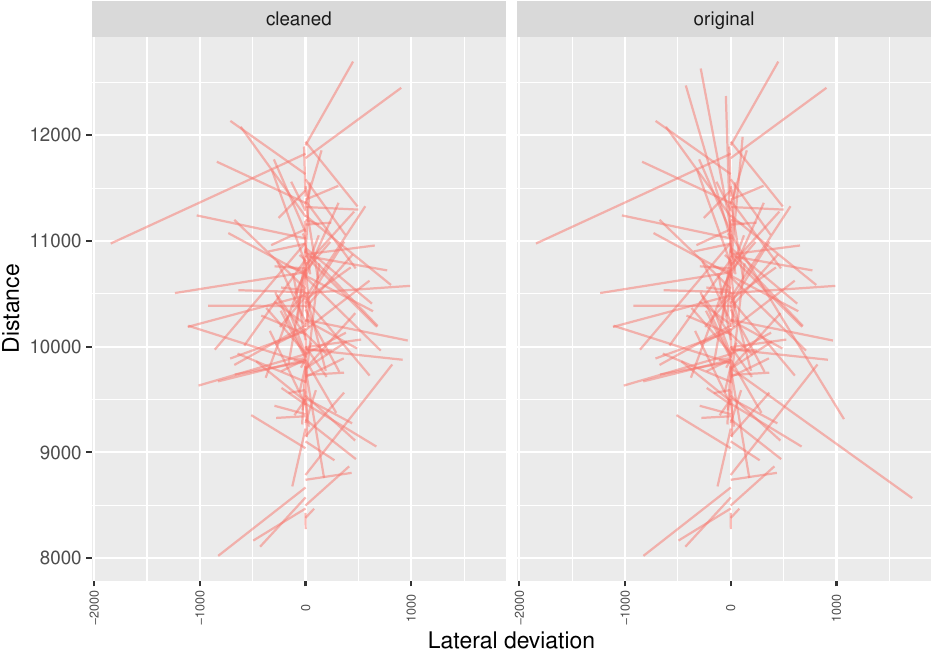} \caption{\label{fig:philtee} Each segment of the figures represents a target/destination pair for shots played off the tee. All target points have been rotated so as to appear on the y-axis. The ``original'' panel shows the inferred data for Phil Mickelson before removing outliers. The ``cleaned'' panel shows the data for Phil Mickelson after removal of outliers. Numbers are in meters.}\label{fig:unnamed-chunk-7}
\end{figure*}

\begin{figure*}[!htbp]
\includegraphics{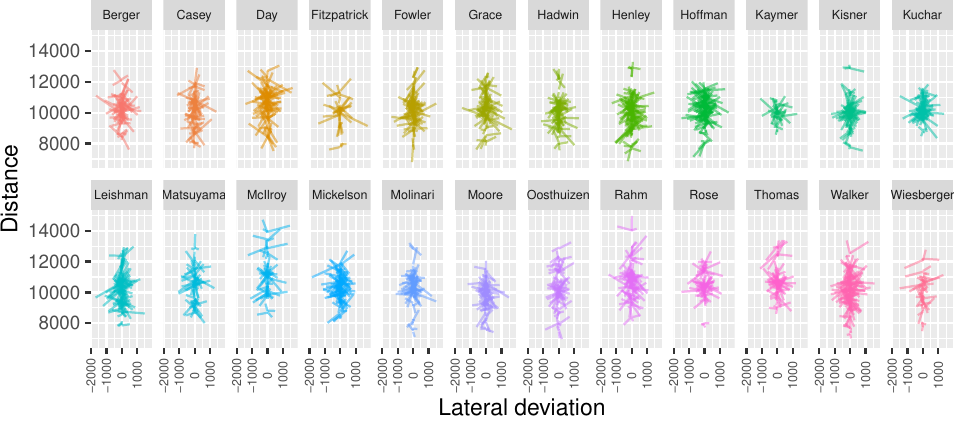} \caption{\label{fig:someplayer1tee} Each segment of the figure represents a target/destination pair for shots played off the tee after removal of outliers. All target points have been rotated so as to appear on the y-axis.  Numbers are in meters.}\label{fig:unnamed-chunk-8}
\end{figure*}

\begin{figure*}[!htbp]
\includegraphics{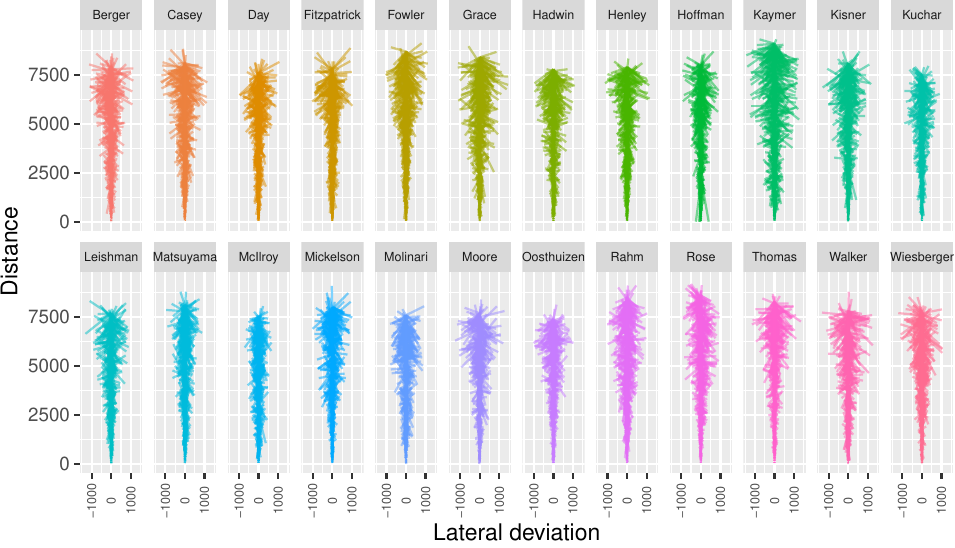} \caption{\label{fig:bootstraprough1} Each segment of the figure represents a target/destination pair for shots played from the rough after bootstrapping. All target points have been rotated so as to appear on the y-axis.  Numbers are in meters.}\label{fig:unnamed-chunk-10}
\end{figure*}

\begin{figure*}[!htbp]
\includegraphics{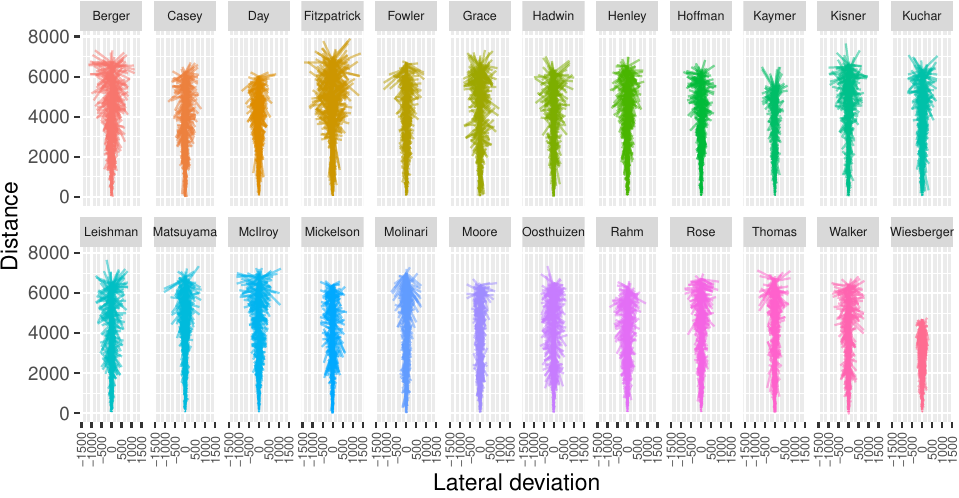} \caption{\label{fig:bootstrapbunker1} Each segment of the figure represents a target/destination pair for shots played from the bunker after bootstrapping. All target points have been rotated so as to appear on the y-axis.  Numbers are in meters.}\label{fig:unnamed-chunk-11}
\end{figure*}

\begin{figure*}[!htbp]
\includegraphics{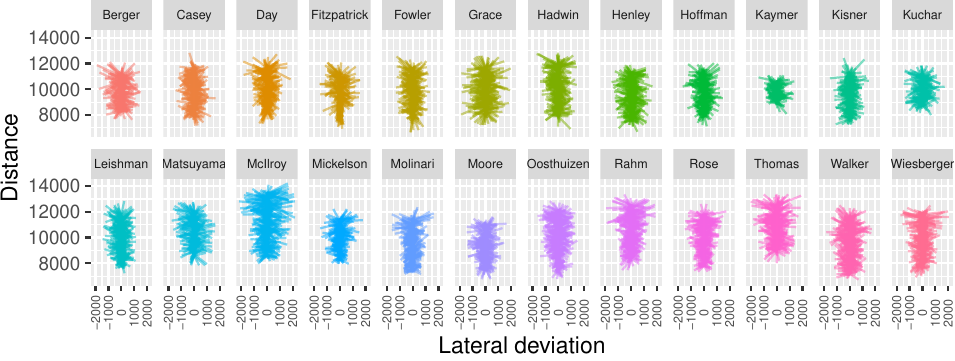} \caption{\label{fig:bootstraptee1} Each segment of the figure represents a target/destination pair for shots played off the tee after bootstrapping. All target points have been rotated so as to appear on the y-axis.  Numbers are in meters.}\label{fig:unnamed-chunk-12}
\end{figure*}
\begin{figure*}
\includegraphics{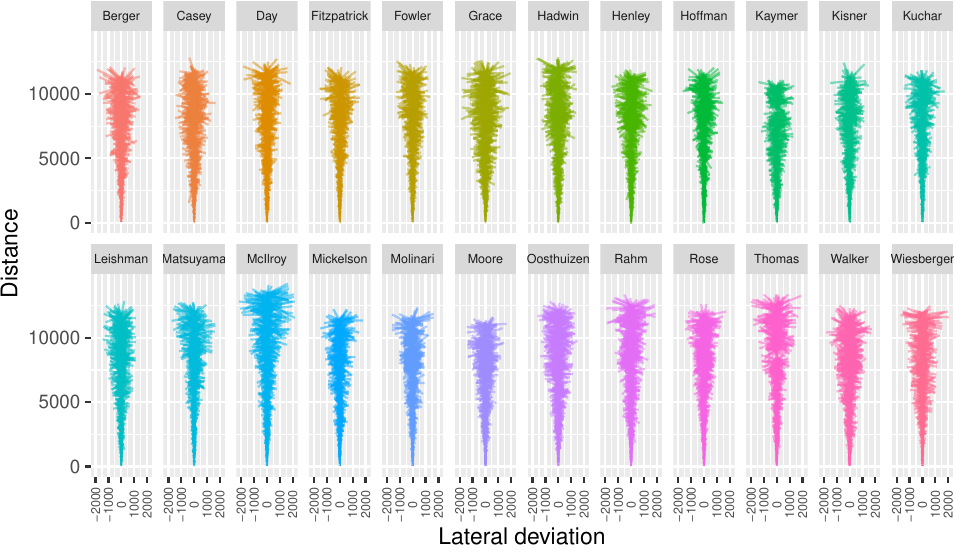} \caption{\label{fig:bootstrapfairwaytee1} Each segment of the figure represents a target/destination pair for shots played off the tee and from the fairway after bootstrapping. All target points have been rotated so as to appear on the y-axis.  Numbers are in meters.}\label{fig:unnamed-chunk-13}
\end{figure*}

\begin{figure}[!htbp]
\includegraphics[width=\textwidth]{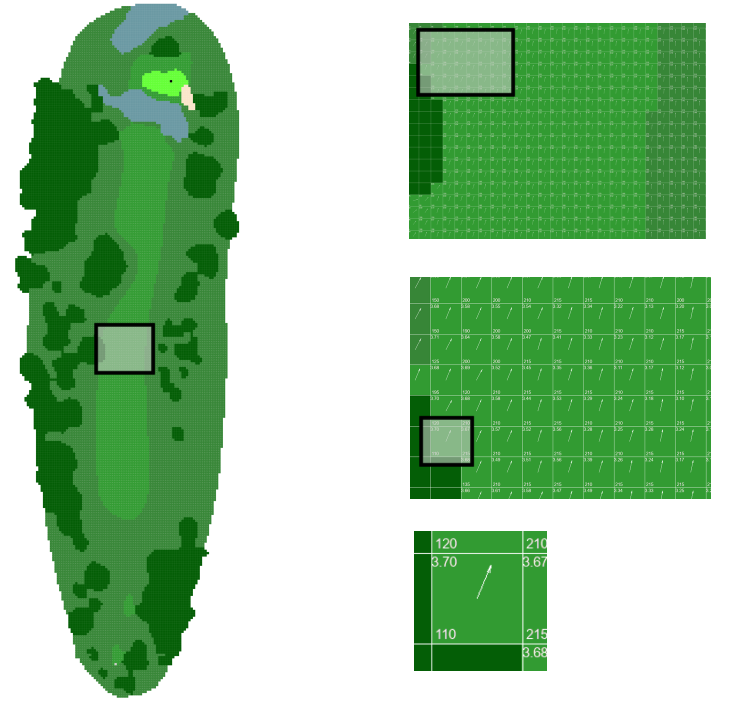}
\caption{\label{fig:application} The figure presents a typical output of the model: for each cell on the golf course, it provides the player with an optimal direction (indicated by an arrow) and a recommended target distance (110 m in this example), along with the expected number of remaining strokes to complete the hole (3.70 in this case). We successively zoom in to reveal finer details at each level.}
\end{figure}

\end{document}